\documentclass[journal]{IEEEtran}
\usepackage{amsmath}    
\usepackage{graphics,epsfig,subfigure}    
\usepackage{verbatim}   
\usepackage{color}      
\usepackage{subfigure}  
\usepackage{hyperref}   
\usepackage{amssymb}
\usepackage{epstopdf}
\usepackage{graphicx}
\newcommand{\prob}{Pr}
\newcommand{\expect}[1]{\mathbb{E}\big\{#1\big\}}

\newcommand{\bv}[1]{{\boldsymbol{#1} }}
\newcommand{\script}[1]{{{\cal{#1} }}}

\begin{document}
\title{Delay Reduction via Lagrange Multipliers in Stochastic Network Optimization}
\author{\large{Longbo Huang, Michael J. Neely}%
\thanks{Longbo Huang (email: longbohu@usc.edu)
and Michael J. Neely (web:  http://www-rcf.usc.edu/$\sim$mjneely)
are with the Department of Electrical
Engineering, University of Southern California, Los Angeles, CA 90089, USA.}%
\thanks{This material is supported in part  by one or more of 
the following: the DARPA IT-MANET program
grant W911NF-07-0028, the NSF grant OCE 0520324, 
the NSF Career grant CCF-0747525.} }
\markboth{Draft}{Huang}
\maketitle

\newtheorem{rem}{Remark}
\newtheorem{fact_def}{\textbf{Fact}}
\newtheorem{coro}{\textbf{Corollary}}
\newtheorem{lemma}{\textbf{Lemma}}
\newtheorem{main}{\textbf{Proposition}}
\newtheorem{thm}{\textbf{Theorem}}
\newtheorem{claim}{\emph{Claim}}
\newtheorem{prop}{Proposition}
\newtheorem{assumption}{\textbf{Assumption}}

\begin{abstract}
In this paper, we consider the problem of reducing network delay in stochastic network utility optimization problems. We start by studying the recently proposed quadratic Lyapunov function based algorithms (QLA). We show that for every stochastic problem, there is a corresponding \emph{deterministic} problem, whose dual optimal solution ``exponentially attracts'' the network backlog process under QLA. In particular, the probability that the backlog vector under QLA deviates from the attractor is exponentially decreasing in their Euclidean distance.  This not only helps to explain how QLA achieves the desired performance but also suggests that one can roughly ``subtract out'' a Lagrange multiplier from the system induced by QLA. We thus develop a family of \emph{Fast Quadratic Lyapunov based Algorithms} (FQLA) that achieve an $[O(1/V), O(\log^2(V))]$ performance-delay tradeoff for problems with a discrete set of action options, and achieve a square-root tradeoff for continuous problems. This is similar to the optimal performance-delay tradeoffs achieved in prior work by Neely (2007) via drift-steering methods, and shows that QLA algorithms can also be used to approach such performance. 

These results highlight the ``network gravity'' role of Lagrange Multipliers in network scheduling. This role can be viewed as the counterpart of the ``shadow price'' role of Lagrange Multipliers in flow regulation for classic flow-based network problems. 

\end{abstract}

\section{Introduction}
In this paper, we consider the problem of reducing network delay in the following general framework of the stochastic network utility optimization problem. We are given a time slotted stochastic network. The network state, such as the network channel condition, is time varying according to some probability law. A network controller performs some action based on the observed network state at every time slot. The chosen action incurs a cost (since cost minimization is mathematically equivalent to utility maximization, below we will use cost and utility interchangeably), but also serves some amount of traffic and possibly generates new traffic for the network.  This traffic causes congestion, and thus leads to backlogs at nodes in the network. The goal of the controller is to minimize its time average cost subject to the constraint that the time average total backlog in the network is finite. 

This setting is very general, and many existing works fall into this category. Further, many techniques have been used to study this problem (see \cite{yichiang_netopt08} for a survey).  In this paper, we focus on algorithms that are built upon quadratic Lyapunov functions (called QLA in the following), e.g., \cite{tassiulas92}, \cite{neelyenergy}, \cite{neelyfairness}, \cite{huangneelypricing}, \cite{rahulneelycognitive}, \cite{neelynowbook}. These QLA algorithms are easy to implement, greedy in nature, and are parameterized by a scalar control variable $V$. It has been shown that when the network state is i.i.d., QLA algorithms can achieve a time average utility that is within $O(1/V)$ to the optimal. Therefore, as $V$ grows large, the time average utility can be pushed arbitrarily close to the optimal. However, such close-to-optimal utility is usually at the expense of large network delay. In fact, in \cite{neelyenergy}, \cite{neelyfairness}, \cite{neelynowbook}, it is shown that an $O(V)$ network delay is incurred when an $O(1/V)$ close-to-optimal utility is achieved. Two recent papers \cite{neelyenergydelay} and \cite{neelysuperfast}, which show that it is possible to achieve within $O(1/V)$ of optimal utility with only $O(\log(V))$ delay, use a more sophisticated algorithm design approach based on exponential Lyapunov functions.  Therefore, it seems that though being simple in implementation, QLA algorithms have undesired delay performance. 


However, we note that the delay results of QLA are usually given in terms of long term upper bounds of the average network backlog e.g., \cite{neelynowbook}. Thus they do not examine the possibility that the actual backlog vector (or its time average) converges to some fixed value. Work in \cite{eryilmaz_qbsc_ton07} considers drift properties towards an ``invariant'' backlog vector, derived in the special case when the problem exhibits a unique optimal Lagrange multiplier. An upper bound on the long term deviation of the actual backlog and the Lagrange multiplier vector is obtained. While this suggests Lagrange multipliers are ``gravitational attractors,'' the bounds in \cite{eryilmaz_qbsc_ton07} do not show that the the actual backlog is very unlikely to deviate significantly from the attractor.

In this paper, we focus on obtaining stronger probability results of the steady state backlog process behavior under QLA. We first show that under QLA, even though the backlog can grow linearly in $V$, it ``typically'' stays close to an ``attractor,'' which is the dual optimal solution of a \emph{deterministic} optimization problem. In particular, the probability that the backlog vector deviates from the attractor is exponentially decreasing in distance, which significantly tightens the attractor analysis in \cite{eryilmaz_qbsc_ton07}. This implies that a large amount of the data is kept in the network simply for maintaining the backlog at the ``right'' level. Therefore, even if we replace these data with some fake data (denoted as \emph{place-holder bits} \cite{neelyrahul_asilomar}), the performance of QLA will not be heavily affected. 
Based on this finding, we propose a family of \emph{Fast Quadratic Lyapunov based Algorithms} (FQLA), which intuitively speaking, can be viewed as subtracting out a Lagrange multiplier from the system induced by QLA.  We show that when the network state is i.i.d., FQLA is able to achieve within  $O(1/V)$ of optimal utility with an $O(\log^2(V))$ delay guarantee for problems with a discrete set of action options, and achieve an $[O(1/V), O(\log^2(V)\sqrt{V})]$ tradeoff for problems with a set of continuous action options. 
The development of FQLA also provides us with additional  insights  into QLA algorithms and the role of Lagrange multipliers in stochastic network optimization. 

The performance of FQLA is closely related to the TOCA algorithm in \cite{neelyenergydelay}, which obtains the same logarithmic and square-root tradeoffs for the  energy-delay problem (up to a $\log(V)$ difference) via drift steering techniques. However, we note that FQLA differs from TOCA in the following: First,  TOCA in \cite{neelyenergydelay} is constructed based on exponential Lyapunov functions; while FQLA uses simpler quadratic Lyapunov functions. Second, FQLA is designed to mimic QLA, thus can be viewed as trying to maintain the dual variable property under QLA; whereas TOCA is designed to ensure the primal constraints are satisfied. Third, FQLA requires an arbitrary small but nonzero fraction of packet droppings, hence can not be applied to problems where packet dropping is not allowed.  

We now summarize the main contributions of this paper in the following:
\begin{itemize}
\item This paper proves that in steady state, the backlog process under QLA is ``exponentially attracted'' to an attractor. This fact also helps to explain how QLA achieves the desired performance.
\item This paper proposes a family of \emph{Fast Quadratic Lyapunov based Algorithms (FQLA)}, which are usually easy to implement, and can achieve an $[O(1/V), O(\log^2(V))]$ performance-delay tradeoff for general stochastic optimization problems with a discrete set of action options as well as a square-root tradeoff for continuous problems. 
\item This paper highlights a new functionality of Lagrange multipliers: the ``network gravity''  in network scheduling. 

\end{itemize}

The paper is organized as follows: In Section \ref{section:notation}, we set up our notations. In Section \ref{section:model}, we state our network model. We then review the QLA algorithm and define the \emph{deterministic problem} in Section \ref{section:qla_rism}. In Section \ref{section:qla_multicon}, we show that the backlog process under QLA always stays close to an attractor. In Section \ref{section:fqla}, we propose the FQLA algorithm. Section \ref{section:secondspecial} considers single queue network problems and provides both \emph{deterministic} and probabilistic bounds on the backlog size. Section \ref{section:numerical} provides simulation results. 
We discuss the ``gravity'' role of Lagrange multipliers and relate QLA to the randomized incremental subgradient method (RISM)  \cite{bertsekasoptbook} in Section \ref{section:LMfunctionality}. 

\section{Notations}\label{section:notation}
\begin{itemize}
\item $\mathbb{R}$:  the set of real numbers
\item $\mathbb{R}_+$ (or $\mathbb{R}_-$): the set of nonnegative (or non-positive) real numbers
\item $\mathbb{R}^n$ (or $\mathbb{R}^n_+$): the set of $n$ dimensional \emph{column} vectors, with each element being in $\mathbb{R}$ (or $\mathbb{R}_+$) 
\item \textbf{bold} symbols $\bv{x}$ and $\bv{x}^T$: \emph{column} vector and its transpose
\item $\bv{x}\succeq\bv{y}$: vector $\bv{x}$ is entrywise no less than vector $\bv{y}$
\item $\bv{0}$: column vector with all elements being $0$
\end{itemize}

$\vspace{-.3in}$
\section{System Model}\label{section:model}
In this section, we specify the general network model we use. We consider a network controller that operates a network with the goal of minimizing the time average cost, subject to the queue stability constraint. The network is assumed to operate in slotted time, i.e., $t\in\{0,1,2,...\}$. We assume there are $r\geq1$ queues in the network. 

$\vspace{-.22in}$
\subsection{Network State}

We assume there are a total of $M$ different random network states, 
and define $\script{S} = \{s_1, s_2, \ldots, s_M\}$ as the set of 
possible states.  Each particular state $s_i$ indicates the current
network parameters, such as a vector of channel conditions for
each link, or a collection of other relevant information about
the current network channels and arrivals.  Let $S(t)$ denote the network 
state at time $t$. We assume that $S(t)$ is i.i.d. every time slot, and let $p_{s_i}$ denote its 
probability of being in state $s_i$, i.e., $p_{s_i}=\prob\{S(t)=s_i\}$. We assume
the network controller can observe $S(t)$ at the beginning of every
slot $t$, but the $p_{s_i}$ probabilities 
are not necessarily known.

$\vspace{-.22in}$
\subsection{The Cost, Traffic and Service}\label{subsection:costtrafficservice}
At each time $t$, after observing $S(t)=s_i$, the controller chooses an action $x(t)$ from a set $\script{X}^{(s_i)}$, i.e., $x(t)= x^{(s_i)}$ for some $x^{(s_i)}\in\script{X}^{(s_i)}$. The set $\script{X}^{(s_i)}$ is called the feasible action set for network state $s_i$ and is assumed to be time-invariant and compact for all $s_i\in\script{S}$.  The cost, traffic and service generated by the chosen action $x(t)=x^{(s_i)}$ are as follows:
\begin{enumerate}
\item[(a)] The chosen action has an associated cost given by the cost function $f(t)=f(s_i, x^{(s_i)}): \script{X}^{(s_i)}\mapsto \mathbb{R}_+$ (or $\script{X}^{(s_i)}\mapsto\mathbb{R}_-$ in the case of reward maximization problems);

\item[(b)] The amount of traffic generated by the action to queue $j$ is determined by the traffic function $A_j(t)=g_{j}(s_i, x^{(s_i)}): \script{X}^{(s_i)}\mapsto \mathbb{R}_{+}$, in units of packets; 

\item[(c)] The amount of service allocated to queue $j$ is given by the rate function $\mu_j(t)=b_{j}(s_i, x^{(s_i)}): \script{X}^{(s_i)}\mapsto \mathbb{R}_{+}$, in units of packets;

 \end{enumerate}
Note that $A_j(t)$ includes both the exogenous arrivals from outside the network to queue $j$, and the endogenous arrivals from other queues, i.e., the transmitted packets from other queues, to queue $j$ (See Section \ref{section:queuenotation} and \ref{subsection:toyexample} for further explanations). We assume the functions $f(s_i, \cdot)$, $g_{j}(s_i, \cdot)$ and $b_{j}(s_i, \cdot)$ are time-invariant, their magnitudes are uniformly upper bounded by some constant $\delta_{max}\in(0,\infty)$ for all $s_i$, $j$, and they are known to the network operator. We also assume that there exists a set of actions $\{x^{(s_i)k}\}_{i=1,..., M}^{k=1,..., r+2}$ with $x^{(s_i)k}\in\script{X}^{(s_i)}$ such that $\sum_{s_i}p_{s_i}\big\{\sum_k\vartheta^{(s_i)}_k[g_{j}(s_i, x^{(s_i)k})-b_{j}(s_i, x^{(s_i)k})]\big\}\leq -\epsilon$ for some $\epsilon>0$ for all $j$, with $\sum_j\vartheta^{(s_i)}_k=1$  and $\vartheta^{(s_i)}_k\geq0$ for all $s_i$ and $k$. That is, the
constraints are feasible with $\epsilon$ slackness. Thus, there exists a stationary randomized policy that stabilizes all queues (where $\vartheta^{(s_i)}_k$ represents the probability of choosing action $x^{(s_i)k}$ when $S(t)=s_i$).  In the following, we use:
\begin{eqnarray}
\bv{A}(t)&=&(A_1(t), A_2(t), ..., A_r(t))^{T},\label{eq:arrivalvector}\\
\bv{\mu}(t)&=&(\mu_1(t), \mu_2(t), ..., \mu_r(t))^{T},\label{eq:servicevector}
\end{eqnarray}
to denote the arrival and service vectors at time $t$. It is easy to see from above that if we define:
\begin{eqnarray}
B=\sqrt{r}\delta_{max},\label{eq:Bdef}
\end{eqnarray}
then $\|\bv{A}(t)-\bv{\mu}(t)\|\leq B$ for all $t$. 

$\vspace{-.25in}$
\subsection{Queueing, Average Cost and the Stochastic Problem}\label{section:queuenotation}
Let $\bv{U}(t)=(U_1(t), ..., U_r(t))^T\in\mathbb{R}^r_{+}$, $t=0, 1, 2, ...$ be the queue backlog vector  process of the network, in units of packets. We assume the following queueing dynamics: 
\begin{eqnarray}
U_j(t+1)=\max\big[U_j(t)-\mu_j(t), 0\big]+A_j(t)\quad\forall j,\label{eq:queuedynamic}
\end{eqnarray}
and $\bv{U}(0)=\bv{0}$. Note that by using (\ref{eq:queuedynamic}), we assume that when a queue does not have enough packets to send, null packets are transmitted.  In this paper, we adopt the following notion of queue stability:
\begin{eqnarray}
\expect{\sum_{j=1}^rU_j}\triangleq
\limsup_{t\rightarrow\infty}\frac{1}{t}\sum_{\tau=0}^{t-1}\sum_{j=1}^{r}\expect{U_j(\tau)}<\infty.\label{eq:queuestable}
\end{eqnarray}
We also use $f^{\pi}_{av}$ to denote the time average cost induced by an action-seeking policy $\pi$, defined as:
\begin{eqnarray}
f^{\pi}_{av}\triangleq
\limsup_{t\rightarrow\infty}\frac{1}{t}\sum_{\tau=0}^{t-1}\expect{f^{\pi}(\tau)},\label{eq:timeavcost}
\end{eqnarray}
where $f^{\pi}(\tau)$ is the cost incurred at time $\tau$ by policy $\pi$. We call an action-seeking  policy under which (\ref{eq:queuestable}) holds a \emph{stable} policy, and use $f_{av}^*$ to denote the optimal time average cost over all stable policies. 
Every slot, the network controller observes the current network state and chooses a control action, with the goal of minimizing time average cost subject to network stability. This goal can be mathematically stated as:
\begin{eqnarray*}
\min: \,\,\, f_{av}, \quad s.t.\,\,\,  (\ref{eq:queuestable}).
\end{eqnarray*}
In the rest of the paper, we will refer to this problem as \emph{the stochastic problem}. This stochastic problem framework can be used to model many network utility problems, such as the energy minimization problem \cite{neelyenergy} and the access point pricing problem \cite{huangneelypricing}.
We note that a similar network model with stochastic penalties is treated in \cite{stolyar06primaldual} using a fluid model and a primal-dual approach that achieves optimality in a limiting sense. The framework is also treated in \cite{neelynowbook} using a quadratic Lyapunov based algorithm (QLA) that provides an explicit $[O(1/V), O(V)]$ performance-delay tradeoff when the network state is i.i.d..

$\vspace{-.2in}$
\subsection{An Example of the Model}\label{subsection:toyexample}
\begin{figure}[cht]
\centering
\includegraphics[height=0.5in, width=2.5in]{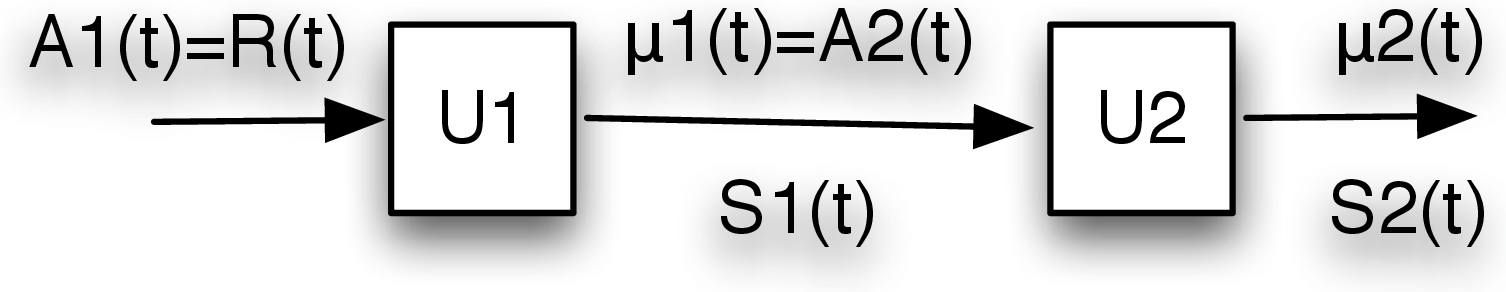}
\caption{A 2-queue system}\label{fig:model_2q}
\end{figure}

Here we provide an example to illustrate our model. Consider the $2$-queue network in Fig.\ref{fig:model_2q}. Every slot, the network operator makes a decision on whether or not to allocate one unit power to serve packets at each queue, so as to support all arriving traffic, i.e., maintain queue stability, with minimum energy expenditure. Every slot, the number of arrival packets $R(t)$, is i.i.d., being either $2$ or $0$ with probabilities $5/8$ and $3/8$ respectively. The channel states $S_1(t),S_2(t)$ are also i.i.d. being either ``G=good'' or ``B=bad'' with equal probabilities. One unit of power can serve $2$ packets in a good channel but can only serve one in a bad channel. Both channels can be activated simultaneously without affecting each other. 

\begin{figure}[cht]
\centering
\includegraphics[height=1in, width=3in]{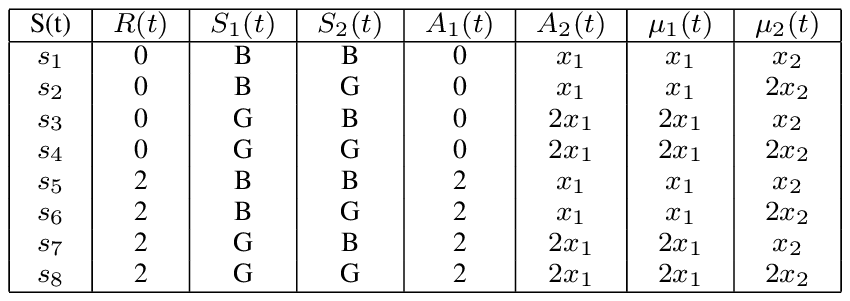}
\caption{Network state, Traffic and Rate functions}\label{fig:exampletable}
\end{figure}

In this case, a network state $S(t)$ is a $(R(t), S_1(t), S_2(t))$ tuple and $S(t)$ is i.i.d.. There are eight possible network states. At each state $s_i$, the action $x^{(s_i)}$ is a pair $(x_1, x_2)$, with $x_i$ being the amount of energy spent at queue $i$, and $(x_1, x_2)\in\script{X}^{(s_i)}=\{0/1,0/1\}$. The cost function is always $f(s_i, x^{(s_i)})=x_1+x_2$ for all $s_i$. The network states, the traffic functions and service rate functions are summarized in Fig. \ref{fig:exampletable}. Note here $A_1(t)=R(t)$ is part of $S(t)$ and thus is independent of $x^{(s_i)}$; while $A_2(t)=\mu_1(t)$ hence depends on $x^{(s_i)}$. Also note that $A_2(t)$ equals $\mu_1(t)$ instead of $\min[\mu_1(t), U_1(t)]$ due to our idle fill assumption in Section \ref{section:queuenotation}.

\vspace{-.2in}
\section{QLA and the Deterministic Problem}\label{section:qla_rism}
In this section, we first review the quadratic Lyapunov functions based algorithms (the QLA algorithm) \cite{neelynowbook} for solving the stochastic problem. Then we define the \emph{deterministic problem} and its dual. We then describe the ordinary subgradient method (OSM) that can be used to solve the dual.  The dual problem and OSM will also be used later for our analysis of the steady state backlog behavior under QLA.

\vspace{-0.1in}
\subsection{The QLA algorithm}
To solve the stochastic problem using QLA, we first define a quadratic Lyapunov function $L(\bv{U}(t))=\frac{1}{2}\sum_{j=1}^{r}U_j^2(t)$. We then define the one-slot conditional Lyapunov drift: $\Delta(\bv{U}(t))=\expect{L(\bv{U}(t+1))-L(\bv{U}(t))\left |\right.\bv{U}(t)}$. From  (\ref{eq:queuedynamic}), we obtain the following drift expression:
\begin{eqnarray*}
\Delta (\bv{U}(t))\leq C - \expect{\sum_{j=1}^{r}U_j(t)\big[\mu_j(t)-A_j(t)\big]\left |\right.\bv{U}(t)},
\end{eqnarray*}
where $C=r\delta_{max}^2$. Now add to both sides the term $V\expect{f(t)\left |\right. \bv{U}(t)}$, where $V\geq1$ is a scalar control variable,  we obtain:
\begin{eqnarray}
\Delta (\bv{U}(t))+V\expect{f(t)\left|\right. \bv{U}(t)}\leq C- \mathbb{E}\bigg\{-Vf(t)  \label{eq:drift}\\
\qquad+\sum_{j=1}^{r}U_j(t)\big[\mu_j(t)-A_j(t)\big]\left|\right.\bv{U}(t)\bigg\}.\nonumber
\end{eqnarray}

The QLA algorithm is then obtained by choosing an action $x$ at every time slot $t$ to \emph{minimize the right hand side of (\ref{eq:drift})} given $\bv{U}(t)$. Specifically, the QLA algorithm works as follows: 


\underline{\emph{QLA:}} At every time slot $t$, observe the current network state $S(t)$ and the backlog $\bv{U}(t)$. If $S(t)=s_i$, choose $x^{(s_i)}\in\script{X}^{(s_i)}$ that solves the following: 
\begin{eqnarray}
\max && -Vf(s_i, x)+\sum_{j=1}^{r}U_j(t)\big[b_j(s_i, x)-g_j(s_i, x)\big]\label{eq:QLAeq}\\
s.t. && x\in\script{X}^{(s_i)}.\nonumber
\end{eqnarray}
Depending on the problem structure, (\ref{eq:QLAeq}) can usually be decomposed into separate parts that are easier to solve, e.g., \cite{neelyenergy}, \cite{huangneelypricing}. 
Also, it can be shown, as in \cite{neelynowbook} that, 
\begin{eqnarray}
f_{av}^{QLA}=f^*_{av}+O(1/V),\quad \overline{U}^{QLA}=O(V),\label{eq:qla_performance}
\end{eqnarray}
where $f_{av}^{QLA}$ is the average cost under QLA and $\overline{U}^{QLA}$ is the time average network backlog size under QLA. 
 
\subsection{The Deterministic Problem}
Consider \emph{the deterministic problem} as follows: 
\begin{eqnarray}
\min&&\script{F}(\bv{x})\triangleq V\sum_{s_i}p_{s_i}f(s_i, x^{(s_i)})\label{eq:primal}\\
s.t.&&\script{G}_j(\bv{x})\triangleq\sum_{s_i}p_{s_i}g_j(s_i, x^{(s_i)})\nonumber\\
&&\qquad\qquad\qquad\leq \script{B}_j(\bv{x})\triangleq\sum_{s_i}p_{s_i}b_j(s_i, x^{(s_i)})\quad\forall\, j\nonumber\\
&& x^{(s_i)}\in \script{X}^{(s_i)}\quad \forall\, i=1, 2, ..., M,\nonumber
\end{eqnarray}
where $p_{s_i}$ corresponds to the probability of $S(t)=s_i$ and $\bv{x}=(x^{(s_1)}, ..., x^{(s_M)})^T$. The dual problem of (\ref{eq:primal}) can be obtained as follows:
\begin{eqnarray}
\hspace{-.0in}\max && q(\bv{U})\label{eq:dualproblem}\\
\hspace{-.0in} s.t.&& \bv{U}\succeq\bv{0}\nonumber,
\end{eqnarray}
where $q(\bv{U})$ is called the dual function and is defined as:
\begin{eqnarray}
q(\bv{U})=\inf_{x^{(s_i)}\in \script{X}^{(s_i)}}\bigg\{ V\sum_{s_i}p_{s_i}f(s_i, x^{(s_i)})\label{eq:dual}\qquad\qquad\qquad\\
+\sum_jU_j\big[ \sum_{s_i}p_{s_i}g_j(s_i, x^{(s_i)})- \sum_{s_i}p_{s_i}b_j(s_i, x^{(s_i)})\big]\bigg\}.\nonumber
\end{eqnarray}
By rearranging the terms, we note that $q(\bv{U})$ can also be written in the following separable form, which is more useful for our later analysis. 
\begin{eqnarray}
 q(\bv{U})=\inf_{x^{(s_i)}\in \script{X}^{(s_i)}}\sum_{s_i}p_{s_i}\bigg\{Vf(s_i, x^{(s_i)})\label{eq:dual_separable}\qquad\qquad\qquad\\
+\sum_jU_j\big[g_j(s_i, x^{(s_i)})- b_j(s_i, x^{(s_i)})\big]\bigg\}.\nonumber
\end{eqnarray}

Here $\bv{U}=(U_1, ..., U_r)^T$ is the \emph{Lagrange multiplier} of (\ref{eq:primal}). It is well known that $q(\bv{U})$ in (\ref{eq:dual}) is concave in the vector $\bv{U}$, and hence the problem (\ref{eq:dualproblem}) can usually be solved efficiently, particularly when cost functions and rate functions are separable over different different network components. It is also well known that in many situations, the optimal value of (\ref{eq:dualproblem}) is the same as the optimal value of (\ref{eq:primal}) and in this case we say that there is no duality gap \cite{bertsekasoptbook}. 

We note that the deterministic problem (\ref{eq:primal}) is not necessarily convex as the sets $\script{X}^{(s_i)}$ are not necessarily convex, and the functions $f(s_i, \cdot)$, $g_j(s_i, \cdot)$ and $b_j(s_i, \cdot)$ are not necessarily convex. Therefore, there may be a duality gap between the deterministic problem (\ref{eq:primal}) and its dual (\ref{eq:dualproblem}). Furthermore,  solving the deterministic problem (\ref{eq:primal}) may not solve the stochastic problem. This is so since  at every network state, the stochastic problem may require time sharing over more than one action, but the solution to the deterministic problem gives only a fixed operating point per network state.  However, one can show, by using an argument similar to showing the existence of an optimal stationary randomized algorithm in \cite{huangneelypricing}, that the dual problem (\ref{eq:dualproblem}) gives the exact  value of $Vf_{av}^*$, where $f_{av}^*$ is the optimal time average cost,  even if (\ref{eq:primal}) is non-convex. 



Among the many algorithms that can be used to solve (\ref{eq:dualproblem}), the following algorithm is the most common one (for performance see \cite{bertsekasoptbook}), we denote it as the \emph{ordinary subgradient method} (OSM):

\underline{\emph{OSM:}} Initialize $\bv{U}(0)$; at every iteration $t$, observe $\bv{U}(t)$, 
\begin{enumerate}
\item Find $x_{\bv{U}}^{(s_i)}\in\script{X}^{(s_i)}$ for $i\in\{1, ..., M\}$ that achieves the infimum of the right hand side of (\ref{eq:dual}). 
\item Using the $\bv{x}_{\bv{U}}=(x_{\bv{U}}^{(s_1)}, x_{\bv{U}}^{(s_)}, ..., x_{\bv{U}}^{(s_M)})^T$ found, update: 
\begin{eqnarray}
\hspace{-.3in}U_j(t+1)=\max\bigg[U_j(t)-\alpha^t\sum_{s_i}p_{s_i}\big[b_j(s_i, x_{\bv{U}}^{(s_i)})\label{eq:osmupdate}\qquad\\
-g_j(s_i, x_{\bv{U}}^{(s_i)})\big],0\bigg].\nonumber
\end{eqnarray}
\end{enumerate}
We use $x^{(s_i)}_\bv{U}$ to highlight its dependency on $\bv{U}(t)$. The term $\alpha^t>0$ is called the \emph{step size} at iteration $t$. In the following, we will always assume $\alpha^t=1$ when referring to OSM. Note that if there is only one network state, QLA and OSM will choose the same action given the same $\bv{U}$, and they differ only by (\ref{eq:queuedynamic}) and (\ref{eq:osmupdate}). The term $\bv{G}_{\bv{U}}=(G_{\bv{U},1}, G_{\bv{U},2}, ..., G_{\bv{U},r})^{T}$,  with: 
\begin{eqnarray}
G_{\bv{U},j}&=&\script{G}_j(\bv{x}_{\bv{U}})-\script{B}_j(\bv{x}_{\bv{U}})\label{eq:subgradient_def}\\
&=&\sum_{s_i}p_{s_i}\big[-b_j(s_i, x_{\bv{U}}^{(s_i)})+ g_j(s_i, x_{\bv{U}}^{(s_i)})\big],\nonumber
\end{eqnarray}
is called the \emph{subgradient} of $q(\bv{U})$ at $\bv{U}(t)$. It is well known that for any other $\hat{\bv{U}}\in\mathbb{R}^r$, we have:
\begin{eqnarray}
(\hat{\bv{U}}-\bv{U}(t))^{T}\bv{G}_{\bv{U}}\geq q(\hat{\bv{U}})-q(\bv{U}(t)).\label{eq:subgradientpro}
\end{eqnarray}
Using $\|\bv{G}_{\bv{U}}\|\leq B$, we note that (\ref{eq:subgradientpro}) also implies:
\begin{eqnarray}
q(\hat{\bv{U}})-q(\bv{U}(t))\leq B\|\hat{\bv{U}}-\bv{U}(t)\|\quad\forall\,\hat{\bv{U}}, \bv{U}\in\mathbb{R}^r\label{eq:dualbddslope}
\end{eqnarray}

We are now ready to study the steady state behavior of $\bv{U}(t)$ under QLA. To simplify notations and highlight the scaling effect of the scalar $V$ in QLA, we use the following notations:
\begin{enumerate}
\item We use $q_0(\bv{U})$ and $\bv{U}_0^*$ to denote the dual objective function and an optimal solution of (\ref{eq:dualproblem}) when $V=1$; and use $q(\bv{U})$ and $\bv{U}_V^*$ (also called the optimal Lagrange multiplier) for their counterparts with general $V\geq1$;
\item We use $x^{(s_i)}_{\bv{U}}$ to denote an action chosen by QLA for a given $\bv{U}(t)$ and  $S(t)=s_i$; and use $\bv{x}_{\bv{U}}=(x^{(s_1)}_{\bv{U}}, ..., x^{(s_M)}_{\bv{U}})^T$ to denote a solution chosen by OSM for a given $\bv{U}(t)$.
\end{enumerate}

To simplify analysis, we assume the following throughout:

\begin{assumption}\label{assumption:nonzeroopt}
$\bv{U}^*_V=(U^*_{V1}, ..., U^*_{Vr})^T$ is unique for all $V\geq1$. 
\end{assumption}

Note that Assumption \ref{assumption:nonzeroopt} is not very restrictive. In fact, it holds in many network utility optimization problems, e.g., \cite{eryilmaz_qbsc_ton07}. In many cases, we also have $\bv{U}^*_V\neq\bv{0}$. Moreover, for the assumption to hold for all $V\geq1$, it suffices to have just $\bv{U}^*_0$ being unique. This is shown in the following lemma regarding the scaling effect of the parameter $V$ on the optimal Lagrange multiplier. 
 
\begin{lemma}\label{lemma:scalingeffectV}
$\bv{U}^*_V=V\bv{U}^*_0.$
\end{lemma}
\begin{proof} 
From (\ref{eq:dual_separable}) we see that:
\begin{eqnarray*}
q(\bv{U})/V=\inf_{x^{(s_i)}\in \script{X}^{(s_i)}}\sum_{s_i}p_{s_i}\bigg\{f(s_i, x^{(s_i)})\qquad\qquad\qquad\\
+\sum_j\hat{U}_j\big[g_{j}(s_i, x^{(s_i)})- b_{j}(s_i, x^{(s_i)})\big]\bigg\},
\end{eqnarray*}
where $\hat{U}_j=\frac{U_j}{V}$. However, the right hand side is exactly $q_0(\hat{\bv{U}})$, and thus  is maximized at $\hat{\bv{U}}=\bv{U}_0^*$. Hence $q(\bv{U})$ is maximized at $V\bv{U}_0^*$.
\end{proof}

\section{Backlog vector behavior under QLA}\label{section:qla_multicon}
In this section we study the backlog vector behavior under QLA of the stochastic problem. We first look at the case when $q_0(\bv{U})$ is ``locally polyhedral.'' We show that $\bv{U}$ is mostly within $O(\log(V))$ distance from $\bv{U}^*_V$ in this case, even when $S(t)$ evolves according to a more general time homogeneous Markovian process. 
We then consider the case when $q_0(\bv{U})$ is ``locally smooth'', and show that $\bv{U}$ is mostly within $O(\sqrt{V}\log(V))$ distance from $\bv{U}^*_V$. As we will see, these two results also explain how QLA functions. 

\subsection{When $q_0()$ is ``locally polyhedral''}\label{section:localpolyhedral}
In this section, we study the backlog vector behavior under QLA for the case where $q_0(\bv{U})$ is \emph{locally polyhedral} with parameters $\epsilon, L$, i.e., there exist $\epsilon, L>0$, , such that for all $\bv{U}\succeq \bv{0}$ with $\|\bv{U}-\bv{U}^*_0\|<\epsilon$, the dual function $q_0(\bv{U})$ satisfies:
\begin{eqnarray}
q_0(\bv{U}^*_0)\geq q_0(\bv{U})+L\|\bv{U}^*_0-\bv{U}\|\label{eq:multiple_dual_condition}
\end{eqnarray} 
We will show that in this case, even if $S(t)$ is a general time homogeneous Markovian process, the backlog vector will mostly be within $O(\log(V))$ distance to $\bv{U}^*_V$. Hence the same is also true when $S(t)$ is i.i.d..

To start, we assume for this subsection that $S(t)$ evolves according to a time homogeneous markovian process. Now we define the following notations. Given $t_0$, define $\script{T}_{s_i}(t_0, k)$ to be the set of slots at which $S(\tau)=s_i$ for $\tau\in[t_0, t_0+k-1]$. For a given $\nu>0$, define the \emph{convergent interval} $T_{\nu}$ \cite{neelypowerjsac} for the $S(t)$ process to be the smallest number of slots such that for any $t_0$, regardless of past history, we have:
\begin{eqnarray}
\sum_{i=1}^M \bigg| p_{s_i}-\frac{\expect{||\script{T}_{s_i}(t_0, T_{\nu})|| \left |\right. \script{H}(t_0)}}{T_{\nu}} \bigg|\leq \nu,
\end{eqnarray}
here $||\script{T}_{s_i}(t_0, T_{\nu})||$ is the cardinality of $\script{T}_{s_i}(t_0, T_{\nu})$, and $\script{H}(t_0)=\{S(\tau)\}_{\tau=0}^{t_0-1}$ denotes the network state history up to time $t_0$. For any $\nu>0$, such a $T_{\nu}$ must exist for any stationary ergodic processes with finite state space, thus $T_{\nu}$ exists for $S(t)$ in particular. When $S(t)$ is i.i.d. every slot, we have $T_{\nu}=1$ for all $\nu\geq0$, as $\expect{||\script{T}_{s_i}(t_0, 1)|| \left |\right. \script{H}(t_0)}=p_{s_i}$. Intuitively, $T_{\nu}$ represents the time needed for the process to reach its ``near'' steady state. 

The following theorem summarizes the main results. Recall that $B$ is defined in (\ref{eq:Bdef}) as the upper bound of the magnitude change of $\bv{U}$ in a slot.
\begin{thm}\label{theorem:prob_multicon_polyhedral}
If $q_0(\bv{U})$ is locally polyhedral with constants $\epsilon, L>0$, independent of $V$, then under QLA, 
\begin{enumerate}
\item[(a)] There exist constants $\nu>0$, $D\geq\eta>0$, all independent of $V$, such that $D=D(\nu), \eta=\eta(\nu)$, and whenever $\|\bv{U}(t)-\bv{U}_V^*\|\geq D$, we have:
\begin{eqnarray}
\hspace{-.2in}\expect{\|\bv{U}(t+T_{\nu})-\bv{U}^*_V\|\left|\right.\bv{U}(t)}\leq \|\bv{U}(t)-\bv{U}_V^*\|-\eta.\label{eq:qla_exp_drift}
\end{eqnarray}
In particular, the constants $\nu$, $D$ and $\eta$ that satisfy (\ref{eq:qla_exp_drift}) can be chosen as follows: Choose $\nu$ as any constant such that $0<\nu<L/B$. Then choose $\eta$ as any value such that $0<\eta< T_{\nu}(L-B\nu)$. Finally, choose $D$ as: \footnote{It can be seen from (\ref{eq:dualbddslope}) that $B\geq L$. Thus $T_{\nu}B>\eta$. }
\begin{eqnarray}
D= \max\bigg[\frac{(T^2_{\nu}+T_{\nu})B^2-\eta^2}{2T_{\nu}(L-\frac{\eta}{T_{\nu}}-B\nu)}, \eta\bigg].\label{eq:qla_Dvalue}
\end{eqnarray}
\item[(b)] For given constants $\nu, D,\eta$ in (a), there exist some constants $c^*, \beta^*>0$, independent of $V$, such that:
\begin{eqnarray}
\script{P}(D, m)\leq c^*e^{-\beta^*m},\label{eq:pm_ineq}
\end{eqnarray}
where $\script{P}(D, m)$ is defined as:
\begin{eqnarray}
\hspace{-.35in}\script{P}(D, m)\triangleq\lim\sup_{t\rightarrow\infty}\frac{1}{t}\sum_{\tau=0}^{t-1}\prob\{\|\bv{U}(\tau)-\bv{U}_V^*\|>D+m\}. \label{eq:pm_def}
\end{eqnarray}
\end{enumerate}
\end{thm}

Note that if $m=\frac{\log(V)}{\beta^*}$, by (\ref{eq:pm_ineq}) we have $\script{P}(D, m)\leq\frac{c^*}{V}$. 
Also if a steady state distribution of $\|\bv{U}(t)-\bv{U}_V^*\|$ exists under QLA, i.e., the limit of $\frac{1}{t}\sum_{\tau=0}^{t-1}\prob\{\|\bv{U}(\tau)-\bv{U}_V^*\|>D+m\}$ exists as $t\rightarrow\infty$,  then one can replace $\script{P}(D, m)$ with the steady state probability that $\bv{U}(t)$ deviates from $\bv{U}^*_V$ by an amount of $D+m$, i.e., $\prob\{\|\bv{U}(t)-\bv{U}_V^*\|>D+m\}$. Therefore Theorem \ref{theorem:prob_multicon_polyhedral} can be viewed as showing that when (\ref{eq:multiple_dual_condition}) is satisfied, for a large $V$, the backlog $\bv{U}(t)$ under QLA will mostly be within $O(\log(V))$ distance from $\bv{U}_V^*$. This implies that the average backlog will roughly be $\sum U^*_{Vj}$, which is typically $\Theta (V)$ by Lemma \ref{lemma:scalingeffectV}. However, this fact will also allow us to build FQLA upon QLA to ``subtract out'' roughly $\sum U^*_{Vj}$ data from the network and reduce network delay. Theorem \ref{theorem:prob_multicon_polyhedral} also highlights a deep connection between the steady state behavior of the network backlog process $\bv{U}(t)$ and the structure of the dual function $q_0(\bv{U})$. We note that (\ref{eq:multiple_dual_condition}) is not very restrictive. In fact, if $q_0(\bv{U})$ is polyhedral (e.g., $\script{X}^{(s_i)}$ is finite for all $s_i$), with a unique optimal solution $\bv{U}^*_0\succeq\bv{0}$, then (\ref{eq:multiple_dual_condition}) can be satisfied (see Section \ref{section:numerical} for an example). To prove the theorem, we need the following lemma.

\begin{lemma}\label{lemma:expecteddistance}
For any $\nu>0$, under QLA, we have for all $t$,
\begin{eqnarray}
&&\hspace{-.3in}\expect{\|\bv{U}(t+T_{\nu})-\bv{U}^*_V\|^2\left |\right. \bv{U}(t)}\label{eq:qla_general_expecteddist}\\
&&\qquad\leq \|\bv{U}(t)-\bv{U}^*_V\|^2 + (T^2_{\nu}+T_{\nu})B^2\nonumber\\
&&\qquad-2T_{\nu}\big(q(\bv{U}^*_V)-q(\bv{U}(t))\big)+2T_{\nu}\nu B\|\bv{U}^*_V-\bv{U}(t)\|.\nonumber
\end{eqnarray}
\end{lemma}
\begin{proof}
See Appendix A.
\end{proof}

We now use Lemma \ref{lemma:expecteddistance} to prove Theorem \ref{theorem:prob_multicon_polyhedral}. 


\begin{proof} (Theorem \ref{theorem:prob_multicon_polyhedral})
Part (a): We first show that if (\ref{eq:multiple_dual_condition}) holds for $q_0(\bv{U})$ with $L$, then it also holds for $q(\bv{U})$ with the same $L$. To this end, suppose  (\ref{eq:multiple_dual_condition})  holds for $q_0(\bv{U})$ for all $\bv{U}$ satisfying $\|\bv{U}-\bv{U}^*_0\|<\epsilon$. Then for any $\bv{U}\succeq\bv{0}$ such that $\|\bv{U}-\bv{U}^*_V\|<\epsilon V$, we have $\|\bv{U}/V-\bv{U}^*_0\|<\epsilon$, hence:
\begin{eqnarray*}
q_0(\bv{U}^*_0)\geq q_0(\bv{U}/V)+L\|\bv{U}^*_0-\bv{U}/V\|.
\end{eqnarray*}
Multiplying both sides by $V$, we get:
\begin{eqnarray*}
Vq_0(\bv{U}^*_0)\geq Vq_0(\bv{U}/V)+LV\|\bv{U}^*_0-\bv{U}/V\|.
\end{eqnarray*}
Now using $\bv{U}^*_V=V\bv{U}^*_0$ and $q(\bv{U})=Vq_0(\bv{U}/V)$, we have for all $\|\bv{U}-\bv{U}^*_V\|<\epsilon V$:
 
\begin{eqnarray}
q(\bv{U}^*_V)\geq q(\bv{U})+L\|\bv{U}^*_V-\bv{U}\|.\label{eq:multiple_dual_condition_q}
\end{eqnarray}
Since $q(\bv{U})$ is concave, we see that (\ref{eq:multiple_dual_condition_q}) indeed holds for all $\bv{U}\succeq\bv{0}$. Now for a given $\eta>0$, if: 
 \begin{eqnarray}
(T^2_{\nu}+T_{\nu})B^2-2T_{\nu}\big(q(\bv{U}^*_V)-q(\bv{U}(t))\big)\label{eq:drift_multiple_norm_con}\qquad\qquad\qquad\\
+2T_{\nu}\nu B\|\bv{U}^*_V-\bv{U}(t)\|\leq\eta^2-2\eta\|\bv{U}^*_V-\bv{U}(t)\|,\nonumber
\end{eqnarray}
then by  (\ref{eq:qla_general_expecteddist}), we have:
\begin{eqnarray*}
\expect{\|\bv{U}(t+T_{\nu})-\bv{U}_V^*\|^2\left|\right. \bv{U}(t)}\leq (\|\bv{U}(t)-\bv{U}_V^*\|-\eta)^2,
\end{eqnarray*}
which then by Jensen's inequality implies:
\begin{eqnarray*}
(\expect{\|\bv{U}(t+T_{\nu})-\bv{U}^*_V\| \left|\right. \bv{U}(t)})^2
\leq(\|\bv{U}(t)-\bv{U}^*_V\| -\eta)^2.
\end{eqnarray*}
Thus (\ref{eq:qla_exp_drift}) follows whenever (\ref{eq:drift_multiple_norm_con}) holds and $\|\bv{U}(t)-\bv{U}_V^*\|\geq\eta$. It suffices to choose $D$ and $\eta$ such that $D\geq\eta$ and that (\ref{eq:drift_multiple_norm_con}) holds whenever $\|\bv{U}(t)-\bv{U}_V^*\|\geq D$. Now note that (\ref{eq:drift_multiple_norm_con}) can be rewritten as the following inequalty:
\begin{eqnarray}
q(\bv{U}^*_V)\geq q(\bv{U}(t))+(B\nu+\frac{\eta}{T_{\nu}})\|\bv{U}^*_V-\bv{U}(t)\|+\script{Y}\label{eq:multiple_dual_condition_alter}
\end{eqnarray}
where $\script{Y}=\frac{(T^2_{\nu}+T_{\nu})B^2-\eta^2}{2T_{\nu}}$. Choose any $\nu>0$ independent of $V$ such that $B\nu<L$ and choose $\eta\in(0, T_{\nu}(L-B\nu))$. By (\ref{eq:multiple_dual_condition_q}), if:
\begin{eqnarray}
L\|\bv{U}(t)-\bv{U}_V^*\|\geq (B\nu+\frac{\eta}{T_{\nu}})\|\bv{U}^*_V-\bv{U}(t)\|+\script{Y}\label{eq:multiple_dual_condition_tomeet_last}
\end{eqnarray}
then (\ref{eq:multiple_dual_condition_alter}) holds. Now choose $D$ as defined in (\ref{eq:qla_Dvalue}), we see that if $\|\bv{U}(t)-\bv{U}_V^*\|\geq D$, then (\ref{eq:multiple_dual_condition_tomeet_last}) holds, which implies  (\ref{eq:multiple_dual_condition_alter}), and equivalently (\ref{eq:drift_multiple_norm_con}). We also have $D\geq\eta$, hence (\ref{eq:qla_exp_drift}) holds.

Part (b): Now we show that (\ref{eq:qla_exp_drift}) implies (\ref{eq:pm_ineq}). Choose constants $\nu$, $D$ and $\eta$ that are independent of $V$ in (a). Denote $Y(t)=\|\bv{U}(t)-\bv{U}^*_V\|$, we see then whenever $Y(t)\geq D$, we have $\expect{Y(t+T_{\nu})-Y(t)\left|\right.\bv{U}(t)}\leq -\eta$. It is also easy to see that $|Y(t+T_{\nu})-Y(t)|\leq T_{\nu}B$, as $B$ is defined in (\ref{eq:Bdef}) as the upper bound of the magnitude change of $\bv{U}$ in a slot. Define $\tilde{Y}(t)=\max\big[Y(t)-D, 0\big]$. We see that whenever $\tilde{Y}(t)\geq T_{\nu}B$, we have:
\begin{eqnarray}
\expect{\tilde{Y}(t+T_{\nu})-\tilde{Y}(t)\left|\right. \bv{U}(t)}\label{eq:Ytilde_sameasY}\qquad\qquad\qquad\qquad\qquad\\
=\expect{Y(t+T_{\nu})-Y(t)\left|\right. \bv{U}(t)}\leq -\eta.\nonumber
\end{eqnarray}
Now define a Lyapunov function of $\tilde{Y}(t)$ to be $L(\tilde{Y}(t))=e^{w\tilde{Y}(t)}$ with some $w>0$, and define the $T_{\nu}$-slot conditional \emph{drift} to be: 
\begin{eqnarray}
\Delta_{T_{\nu}}(\tilde{Y}(t))&\triangleq&\expect{L(\tilde{Y}(t+T_{\nu}))-L(\tilde{Y}(t))\left|\right. \bv{U}(t)}\nonumber\\
&=&\expect{e^{w\tilde{Y}(t+T_{\nu})}-e^{w\tilde{Y}(t)}\left|\right.\bv{U}(t)}.\label{eq:drift_exp_tildeY}
\end{eqnarray}
It is shown in Appendix B that by choosing $w=\frac{\eta}{T^2_{\nu}B^2+T_{\nu}B\eta/3}$, we have for all $\tilde{Y}(t)\geq0$:
\begin{eqnarray}
\Delta_{T_{\nu}}(\tilde{Y}(t)) &\leq&e^{2wT_{\nu}B}-\frac{w\eta}{2}e^{w\tilde{Y}(t)}. \label{eq:drift_tildeY_combine}
\end{eqnarray}
Taking expectation on both sides, we have:
\begin{eqnarray}
\expect{e^{w\tilde{Y}(t+T_{\nu})}-e^{w\tilde{Y}(t)}}\leq e^{2wT_{\nu}B}-\frac{w\eta}{2}\expect{e^{w\tilde{Y}(t)}}. \label{eq:driftafterexpection_fortelesum}
\end{eqnarray}
Now summing (\ref{eq:driftafterexpection_fortelesum}) over $t\in\{t_0, t_0+T_{\nu}, ..., t_0+(N-1)T_{\nu}\}$ for some $t_0\in\{0, 1, ..., T_{\nu}-1\}$, we have:
\begin{eqnarray*}
\expect{e^{w\tilde{Y}(t_0+NT_{\nu})}-e^{w\tilde{Y}(t_0)}}\leq Ne^{2wT_{\nu}B}\label{eq:Tslotdrift_sumoverj}\qquad\qquad\qquad\quad\\
-\sum_{j=0}^{N-1}\frac{w\eta}{2}\expect{e^{w\tilde{Y}(t_0+jT_{\nu})}}. \nonumber
\end{eqnarray*}
Rearrange the terms, we have:
\begin{eqnarray*}
\sum_{j=0}^{N-1}\frac{w\eta}{2}\expect{e^{w\tilde{Y}(t_0+jT_{\nu})}}\leq Ne^{2wT_{\nu}B}+\expect{e^{w\tilde{Y}(t_0)}}.
\end{eqnarray*}
Summing the above over $t_0\in\{0,1,..., T_{\nu}-1\}$, we obtain:
\begin{eqnarray*}
\sum_{t=0}^{NT_{\nu}-1}\frac{w\eta}{2}\expect{e^{w\tilde{Y}(t)}}\leq NT_{\nu}e^{2wT_{\nu}B}+\sum_{t_0=0}^{T_{\nu}-1}\expect{e^{w\tilde{Y}(t_0)}}.
\end{eqnarray*}
Dividing both sides with $NT_{\nu}$, we obtain:
\begin{eqnarray}
\frac{1}{NT_{\nu}}\sum_{t=0}^{NT_{\nu}-1}\frac{w\eta}{2}\expect{e^{w\tilde{Y}(t)}} \leq e^{2wT_{\nu}B}\label{eq:initialnotmatter1}\qquad\qquad\qquad\\
+\frac{1}{NT_{\nu}}\sum_{t_0=0}^{T_{\nu}-1}\expect{e^{w\tilde{Y}(t_0)}}.\nonumber
\end{eqnarray}
Taking the limsup as $N$ goes to infinity, we obtain:
\begin{eqnarray}
\lim\sup_{t\rightarrow\infty}\frac{1}{t}\sum_{\tau=0}^{t-1}\frac{w\eta}{2}\expect{e^{w\tilde{Y}(\tau)}} &\leq&e^{2wT_{\nu}B}.\label{eq:initialnotmatter2}
\end{eqnarray}
Using the fact that $\expect{e^{w\tilde{Y}(\tau)}}\geq e^{wm}\prob\{\tilde{Y}(\tau)>m\}$, 
\begin{eqnarray}
\hspace{-.2in}\lim\sup_{t\rightarrow\infty}\frac{1}{t}\sum_{\tau=0}^{t-1}\frac{w\eta}{2}e^{wm}\prob\{\tilde{Y}(\tau)>m\} &\leq&e^{2wT_{\nu}B}.\label{eq:prob_tildeY}
\end{eqnarray}
Plug in $w=\frac{\eta}{T_{\nu}^2B^2+T_{\nu}B\eta/3}$ and use the definition of $\tilde{Y}(t)$:
\begin{eqnarray}
\script{P}(D, m) &\leq&\frac{2e^{2wT_{\nu}B}}{w\eta}e^{-wm}\label{eq:prob_tildeY_final_alter}\\
&\leq &\frac{2(T_{\nu}^2B^2+T_{\nu}B\eta/3)e^{\frac{2\eta}{T_{\nu}B+\eta/3}}}{\eta^2} e^{-\frac{\eta m}{T_{\nu}^2B^2+T_{\nu}B\eta/3}}, \nonumber
\end{eqnarray}
where $\script{P}(D, m)$ is defined in (\ref{eq:pm_def}). Therefore (\ref{eq:pm_ineq}) holds with: 
\begin{eqnarray}
c^*&=&\frac{2(T_{\nu}^2B^2+T_{\nu}B\eta/3)e^{\frac{2\eta}{T_{\nu}B+\eta/3}}}{\eta^2},\nonumber\\
\beta^*&=&\frac{\eta}{T_{\nu}^2B^2+T_{\nu}B\eta/3}.\label{eq:pm_ineq_const}
\end{eqnarray}
It is easy to see that $c^*$ and $\beta^*$ are both independent of $V$. 
\end{proof}

Note from (\ref{eq:initialnotmatter1}) and (\ref{eq:initialnotmatter2}) that Theorem \ref{theorem:prob_multicon_polyhedral} indeed holds for any finite $\bv{U}(0)$. We will later use this fact to prove the performance of FQLA. The following theorem is a special case of Theorem  \ref{theorem:prob_multicon_polyhedral} and gives a more direct illustration of Theorem  \ref{theorem:prob_multicon_polyhedral}. Recall that $\script{P}(D,m)$ is defined in (\ref{eq:pm_def}). Define:
\begin{eqnarray}
\hspace{-.3in}&&\script{P}^{(r)}(D, m)\label{eq:pmr_def}\\
\hspace{-.3in}&&\qquad\quad\triangleq\lim\sup_{t\rightarrow\infty}\frac{1}{t}\sum_{\tau=0}^{t-1}\prob\{\exists\, j, |U_j(\tau)-U_{Vj}^*|>D+m\}.\nonumber
\end{eqnarray}
\begin{thm}\label{thm:prob_multi_con}
If the condition in Theorem  \ref{theorem:prob_multicon_polyhedral} holds and $S(t)$ is i.i.d., then under QLA, for any $c>0$:
\begin{eqnarray}
\script{P}(D_1, cK_1\log(V))&\leq& \frac{c_1^*}{V^c},\label{eq:prob_pm_special}\\
\script{P}^{(r)}(D_1, cK_1\log(V))&\leq& \frac{c_1^*}{V^c}.\label{eq:prob_pmr_special}
\end{eqnarray}
where $D_1=\frac{2B^2}{L}+\frac{L}{4}$, $K_1=\frac{B^2+BL/6}{L/2}$ and $c_1^*=\frac{8(B^2+BL/6)e^{\frac{L}{B+L/6}}}{L^2}$.
\end{thm}
\begin{proof}
First we note that when $S(t)$ is i.i.d., we have $T_{\nu}=1$ for $\nu=0$. Now choose $\nu=0$, $T_{\nu}=1$ and $\eta=L/2$, then we see from (\ref{eq:qla_Dvalue}) that 
\[D=\max\bigg[\frac{2B^2-L^2/4}{L}, \frac{L}{2}\bigg]\leq \frac{2B^2}{L}+\frac{L}{4}.\]  
Now by (\ref{eq:pm_ineq_const}) we see that (\ref{eq:pm_ineq}) holds with $c^*=c_1^*$ and $\beta^*=\frac{L/2}{B^2+BL/6}$. Thus by taking $D_1=\frac{2B^2}{L}+\frac{L}{4}$, we have:
\begin{eqnarray*}
\script{P}(D_1, cK_1\log(V))&\leq& c^*e^{-cK_1\beta^*\log(V)}\\
&=&c_1^*e^{-c\log(V)},
\end{eqnarray*}
where the last step follows since $\beta^*K_1=1$. Thus (\ref{eq:prob_pm_special}) follows. Equation (\ref{eq:prob_pmr_special}) follows from (\ref{eq:prob_pm_special}) by using the fact that for any constant $\zeta$, the events $\script{E}_1=\{\exists\, j, |U_j(\tau)-U_{Vj}^*|>\zeta\}$ and $\script{E}_2=\{\|\bv{U}(\tau)-\bv{U}_{V}^*\|>\zeta\}$ satisfy $\script{E}_1\subset\script{E}_2$. Thus:
$\prob\{\exists\, j, |U_j(\tau)-U_{Vj}^*|>\zeta\}\leq \prob\{ \|\bv{U}(\tau)-\bv{U}_{V}^*\|>\zeta\}.$
\end{proof}

Theorem \ref{thm:prob_multi_con} can be viewed as showing that for a large $V$, the probability for $U_j(t)$ to deviate from the $j^{th}$ component of $\bv{U}^*_{V}$ is exponentially decreasing in the distance. Thus it rarely deviates from $U^*_{Vj}$ by more than $\Theta(\log(V))$ distance. Note that one can similarly prove the following theorem for OSM:
\begin{thm}\label{coro:drift_multiple_norm}
If  the condition in Theorem \ref{theorem:prob_multicon_polyhedral} holds, then there exist positive constants $D=\Theta(1)$ and $\eta=\Theta(1)$, i.e, independent of $V$,  such that, under OSM, if $ \|\bv{U}(t)-\bv{U}^*_V\|\geq D$, 
\begin{eqnarray}
\|\bv{U}(t+1)-\bv{U}_V^*\|\leq\|\bv{U}(t)-\bv{U}_V^*\|-\eta.\label{eq:drift_multiple_norm}
\end{eqnarray}
\end{thm}
\begin{proof}
It is easy to show that under OSM, Lemma \ref{lemma:expecteddistance} holds with $\nu=0$, $T_{\nu}=1$ and without the expectation. Indeed, by (\ref{eq:osmupdate}),  (\ref{eq:subgradient_def}) and Lemma \ref{lemma:queuedynamic_diff} in Appendix A, we have: 
\begin{eqnarray*}
\|\bv{U}(t+1)-\bv{U}^*_V\|^2&\leq& \|\bv{U}(t)-\bv{U}^*_V\|^2 + 2B^2\nonumber\\
&& \qquad\qquad-2 (\bv{U}^*_V-\bv{U}(t))^T\bv{G}_{\bv{U}}
\end{eqnarray*}
Now by (\ref{eq:subgradientpro}) we have: $(\bv{U}^*_V-\bv{U}(t))^T\bv{G}_{\bv{U}}\geq q(\bv{U}^*_V)-q(\bv{U}(t))$. Plug this into the above equation, we obtain:
\begin{eqnarray*}
\|\bv{U}(t+1)-\bv{U}^*_V\|^2&\leq& \|\bv{U}(t)-\bv{U}^*_V\|^2 + 2B^2\nonumber\\
&&\qquad\qquad-2\big(q(\bv{U}^*_V)-q(\bv{U}(t))\big)
\end{eqnarray*}
The theorem then follows by using the same argument as in the proof of Theorem \ref{theorem:prob_multicon_polyhedral}.
\end{proof}

Therefore, when there is a single network state, we see that given (\ref{eq:multiple_dual_condition}), the backlog process converges to a ball of size $\Theta(1)$ around $\bv{U}^*_V$.

\subsection{When $q_0()$ is ``locally smooth'' }\label{section:firstspecial}
In this section, we consider the backlog behavior under QLA, for the case where the dual function $q_0(\bv{U})$ is ``locally smooth'' at $\bv{U}^*_0$. Specifically, we say that the function $q_0(\bv{U})$ is \emph{locally smooth} at $\bv{U}^*_0$ with parameters $\varepsilon, L>0$ if for all $\bv{U}\succeq\bv{0}$ such that $\|\bv{U}-\bv{U}_0^*\|< \varepsilon$, we have: 
\begin{eqnarray}
q_0(\bv{U}^*_0)\geq q_0(\bv{U})+L\|\bv{U}-\bv{U}_0^*\|^2,\label{eq:qtwicediff_condition}
\end{eqnarray}
This condition contains the case when $q_0(\bv{U})$ is twice differentiable with $\nabla q(\bv{U}_0^*)=\bv{0}$ and  $\bv{x}^T\nabla^2 q(\bv{U})\bv{x}\leq -2L\|\bv{x}\|^2$ for any $\bv{U}$ with  $\|\bv{U}^*_0-\bv{U}\|< \varepsilon$. Such a case usually occurs when the sets $\script{X}^{(s_i)}, i=1,...,M$ are convex, thus a ``continuous'' set of actions are available. 
Notice that (\ref{eq:qtwicediff_condition}) is a looser condition than (\ref{eq:multiple_dual_condition}) in the neighborhood of $\bv{U}_0^*$. As we will see, such structural difference of $q_0(\bv{U})$ in the neighborhood of $\bv{U}^*_0$ greatly affects the behavior of backlogs under QLA.

\begin{thm}\label{theorem:qtwicediff}
If $q_0(\bv{U})$ is locally smooth at $\bv{U}^*_0$ with parameters $\varepsilon, L>0$, independent of $V$, 
then under QLA with a sufficiently large $V$, we have:
\begin{enumerate}
\item[(a)] There exists $D=\Theta(\sqrt{V})$ such that  whenever $\|\bv{U}-\bv{U}_V^*\|\geq D$, we have:
\begin{eqnarray}
\hspace{-.3in}\expect{\|\bv{U}(t+1)-\bv{U}^*_V\|\left|\right.\bv{U}(t)}\leq \|\bv{U}(t)-\bv{U}_V^*\|-\frac{1}{\sqrt{V}}.\label{eq:qtwicediff_expdrift}
\end{eqnarray}
\item[(b)] $\script{P}(D, m)\leq c^*e^{-\beta^*m}$, where $\script{P}(D, m)$ is defined in (\ref{eq:pm_def}), $c^*=\Theta(V)$ and $\beta^*=\Theta(1/\sqrt{V})$.
\end{enumerate}
\end{thm}

Theorem \ref{theorem:qtwicediff} can be viewed as showing that, when $q_0(\bv{U})$ is locally smooth at $\bv{U}^*_0$, the backlog vector will mostly be within $O(\sqrt{V}\log(V))$ distance from $\bv{U}^*_V$. This contrasts with Theorem \ref{theorem:prob_multicon_polyhedral}, which shows that the backlog will mostly be within $O(\log(V))$ distance from $\bv{U}^*_V$. Intuitively,  this is due to the fact that under local smoothness, the drift towards $\bv{U}^*_V$ is smaller as $\bv{U}$ gets closer to $\bv{U}^*_V$, hence a $\Theta(\sqrt{V})$ distance is needed to guarantee a drift of size $\Theta(1/\sqrt{V})$; whereas under (\ref{eq:multiple_dual_condition}), any nonzero $\Theta(1)$ deviation from $\bv{U}^*_V$ roughly generates a drift of size $\Theta(1)$ towards $\bv{U}^*_V$, ensuring the backlog stays within $O(\log(V))$ distance from $\bv{U}^*_V$.
To prove Theorem \ref{theorem:qtwicediff}, we need the following corollary of Lemma \ref{lemma:expecteddistance}. 
\begin{coro}\label{coro:expecteddistance}
If $S(t)$ is i.i.d., then under QLA,
\begin{eqnarray*}
\expect{\|\bv{U}(t+1)-\bv{U}^*_V\|^2\left |\right. \bv{U}(t)}\leq \|\bv{U}(t)-\bv{U}^*_V\|^2  + 2B^2\nonumber\\
\qquad-2\big(q(\bv{U}^*_V)-q(\bv{U}(t))\big).\nonumber
\end{eqnarray*}
\end{coro}
\begin{proof}
When $S(t)$ is i.i.d., we have $T_{\nu}=1$ for $\nu=0$.
\end{proof}

\begin{proof} (Theorem \ref{theorem:qtwicediff})
Part (a): We first see that for any $\bv{U}$ with  $\|\bv{U}-\bv{U}_V^*\|< \varepsilon V$, we have $\|\bv{U}/V-\bv{U}_0^*\|< \varepsilon$. Therefore, 
\begin{eqnarray}
q_0(\bv{U}^*_0)\geq q_0(\bv{U}/V)+L\|\bv{U}/V-\bv{U}_0^*\|^2.
\end{eqnarray}
Multiply both sides with $V$, we get:
\begin{eqnarray}
q(\bv{U}^*_V)\geq q(\bv{U})+\frac{L}{V}\|\bv{U}-\bv{U}_V^*\|^2.\label{eq:qcon_twicediff_generalV}
\end{eqnarray}
Similar as in the proof of Theorem \ref{theorem:prob_multicon_polyhedral} and by Corollary \ref{coro:expecteddistance}, we see that for (\ref{eq:qtwicediff_expdrift}) to hold, we need $\|\bv{U}(t)-\bv{U}^*_V\|\geq\frac{1}{\sqrt{V}}$ and:
\begin{eqnarray*}
2B^2-2\big(q(\bv{U}^*_V)-q(\bv{U}(t))\big)\leq \frac{1}{V}-\frac{2}{\sqrt{V}}\|\bv{U}(t)-\bv{U}^*_V\|, 
\end{eqnarray*}
which can be rewritten as: 
\begin{eqnarray}
q(\bv{U}^*_V)\geq q(\bv{U}(t))\big)+\frac{1}{\sqrt{V}}\|\bv{U}(t)-\bv{U}^*_V\|+\frac{2B^2-\frac{1}{V}}{2}.\label{eq:qtwicediff_coro_induced}
\end{eqnarray}
By (\ref{eq:qcon_twicediff_generalV}), we see that for (\ref{eq:qtwicediff_coro_induced}) to hold, we only need:
\begin{eqnarray}
\frac{L}{V}\|\bv{U}-\bv{U}_V^*\|^2\geq \frac{1}{\sqrt{V}} \|\bv{U}-\bv{U}_V^*\| +B^2.\label{eq:qtwicediff_power2inequality}
\end{eqnarray}
It is easy to see that (\ref{eq:qtwicediff_power2inequality}) holds whenever:
\begin{eqnarray*}
\|\bv{U}-\bv{U}_V^*\|\geq \frac{\frac{1}{\sqrt{V}}+\sqrt{\frac{1}{V}+ \frac{4B^2L}{V}} }{2L/V} = \frac{\sqrt{V}+\sqrt{V+ 4B^2LV} }{2L}
\end{eqnarray*}
Denote $D=\frac{\sqrt{V}+\sqrt{V+ 4B^2LV} }{2L}$. We see now when $V$ is large, (\ref{eq:qtwicediff_expdrift}) holds for any $\bv{U}$ with $D\leq\|\bv{U}-\bv{U}_V^*\|<\varepsilon V$. Now since $q(\bv{U})$ is concave, it is easy to show that (\ref{eq:qtwicediff_coro_induced}) holds for all $\|\bv{U}-\bv{U}_V^*\|\geq D$. Hence (\ref{eq:qtwicediff_expdrift}) holds for all $\|\bv{U}-\bv{U}_V^*\|\geq D$, proving Part (a). 

Part (b): By an argument that is similar as in the proof of Theorem \ref{theorem:prob_multicon_polyhedral}, we see that Part (b) follows with:
$\beta^*=\frac{3}{3\sqrt{V}B^2+B}$ and $ c^*=2(VB^2+B\sqrt{V}/3)e^{\frac{6}{3B\sqrt{V}+1}}$.
\end{proof}
Notice in this case we can also prove a similar result as Theorem \ref{coro:drift_multiple_norm} for OSM, with the only difference that $D=\Theta(\sqrt{V})$.


\subsection{Implications of Theorem \ref{theorem:prob_multicon_polyhedral} and \ref{theorem:qtwicediff}}

Consider the following simple problem: an operator operates a single queue and tries to support a Bernoulli arrival, i.e., either $1$ or $0$ packet arrives every slot, with rate $\lambda=0.5$ (the rate may be unknown to the operator) with minimum energy expenditure. The channel is time-invariant. The rate-power curve over the channel is given by: $\mu(t)=\log(1+PW(t))$, where $PW(t)$ is the allocated power at time $t$. Thus to obtain a rate of $\mu(t)$, we need $PW(t)=e^{\mu(t)}-1$. Every time slot, the operator decides how much power to allocate and serves the queue at the corresponding rate, with the goal of minimizing the time average power consumption subject to queue stability. Let $\Phi$ denote the time average energy expenditure incurred by the optimal policy. It is not difficult to see that $\Phi=e^{0.5}-1$. 

Now we look at the deterministic problem:
\begin{eqnarray*}
\min:&& V(e^{\mu}-1)\\
s.t.: && 0.5\leq\mu
\end{eqnarray*}
It is easy to obtain $q(U)=\inf_{\mu}\big\{V(e^{\mu}-1)+U(0.5-\mu)\big\}$. Hence by the KKT conditions \cite{bertsekasoptbook} one obtains that  $U_V^*=Ve^{0.5}$ and the optimal policy is to serve the queue at the constant rate $\mu^*=0.5$. Suppose now QLA is applied to the problem. Then, at every slot $t$, given $U(t)=U$, QLA chooses the power to achieve the rate $\mu(t)$ such that: 
\begin{eqnarray}
\mu(t)\in\arg\min \{V(e^{\mu}-1)+U(0.5-\mu)\}=\log(\frac{U(t)}{V}). \label{eq: example_QLA_action}
\end{eqnarray}
which incurs an instantanous power consumption of $PW(t)=\frac{U(t)}{V}$. Now by Theorem \ref{theorem:qtwicediff}, for most of the time $U(t)\in[U^*_V-\sqrt{V}, U^*_V+\sqrt{V}]$, i.e.,  $U(t)\in[Ve^{0.5}-\sqrt{V}, Ve^{0.5}+\sqrt{V}]$. Hence it is almost always the case that:
\[\log(e^{0.5}-\frac{1}{\sqrt{V}})\leq \mu(t)\leq \log(e^{0.5}+\frac{1}{\sqrt{V}}),\]
which implies:
$0.5-\frac{1}{\sqrt{V}}\leq \mu(t)\leq 0.5+\frac{1}{\sqrt{V}}$.
Thus by a similar argument as in \cite{neelyenergydelay}, one can show that $\overline{PW}\leq\Phi+O(1/V)$, where $\overline{PW}$ is the average power consumption. 

Now consider the case when we can only choose to operate at $\mu\in\{0, \frac{1}{4}, \frac{3}{4}, 1\}$, with the corresponding power consumptions being: $PW\in\{0, e^{\frac{1}{4}}-1,e^{\frac{3}{4}}-1,e-1 \}$. One can similarly obtain $\Phi=\frac{1}{2}(e^{\frac{3}{4}}+e^{\frac{1}{4}})$ and $U^*_V=2V(e^{\frac{3}{4}}-e^{\frac{1}{4}})$. In this case, $\Phi$ is achieved by time sharing the two rates $\{\frac{1}{4}, \frac{3}{4}\}$ with equal portion of time. Now by Theorem \ref{theorem:prob_multicon_polyhedral}, we see that under QLA, $U(t)$ is mostly within $\log(V)$ distance to $U^*_V$. Hence by (\ref{eq: example_QLA_action}), we see that QLA  almost always chooses between the two rates $\{\frac{1}{4}, \frac{3}{4}\}$, and uses them with almost equal frequencies. Hence QLA is also able to achieve $\overline{PW}=\Phi+O(1/V)$ in this case.

The above argument can be generalized to many stochastic network optimization problems. Thus, we see that Theorem \ref{theorem:prob_multicon_polyhedral} and \ref{theorem:qtwicediff} not only provide us with probabilistic deviation bounds of $\bv{U}(t)$ from $\bv{U}^*$, but also help to explain why QLA is able to achieve the desired utility performance:  \emph{under QLA,$\bv{U}(t)$ always stays close to $\bv{U}^*_V$, hence the chosen action is always close to the set of optimal actions}. 

\section{The FQLA Algorithm}\label{section:fqla}
In this section, we propose a family of \emph{Fast Quadratic Lyapunov based Algorithms} (FQLA) for general stochastic network optimization problems. We first provide an  example to illustrate the idea of FQLA. We then describe FQLA with known $\bv{U}_V^*$, called FQLA-Ideal, and study its performance. After that, we describe the more general FQLA without such knowledge, called FQLA-General. For brevity, we only describe FQLA for the case when $q_0(\bv{U})$ is locally polyhedral. FQLA for the other case is briefly discussed in Section \ref{section:fqla_smooth}.

\vspace{-.1in}
\subsection{FQLA: a Single Queue Example}
\begin{figure}[cht]
\centering
\includegraphics[height=1.15in, width=1.5in]{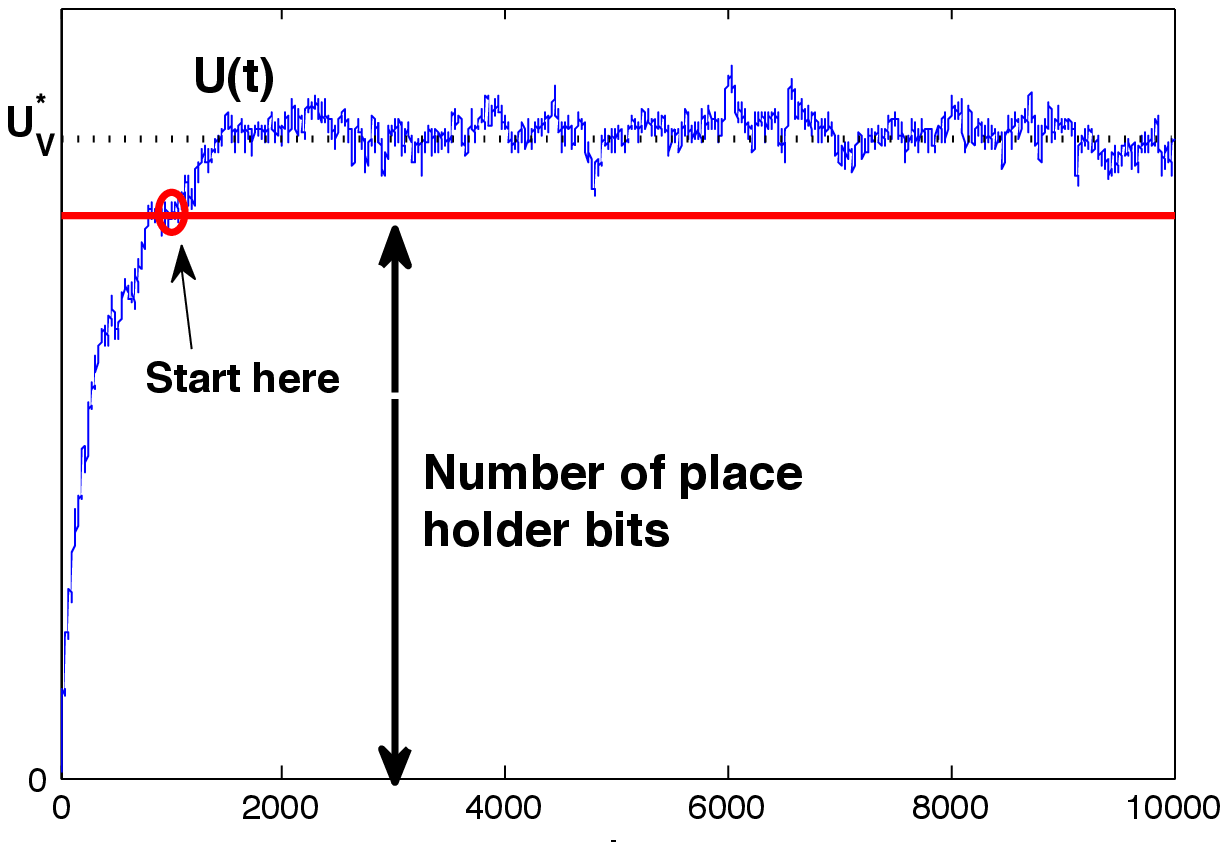}
\includegraphics[height=1.14in, width=1.5in]{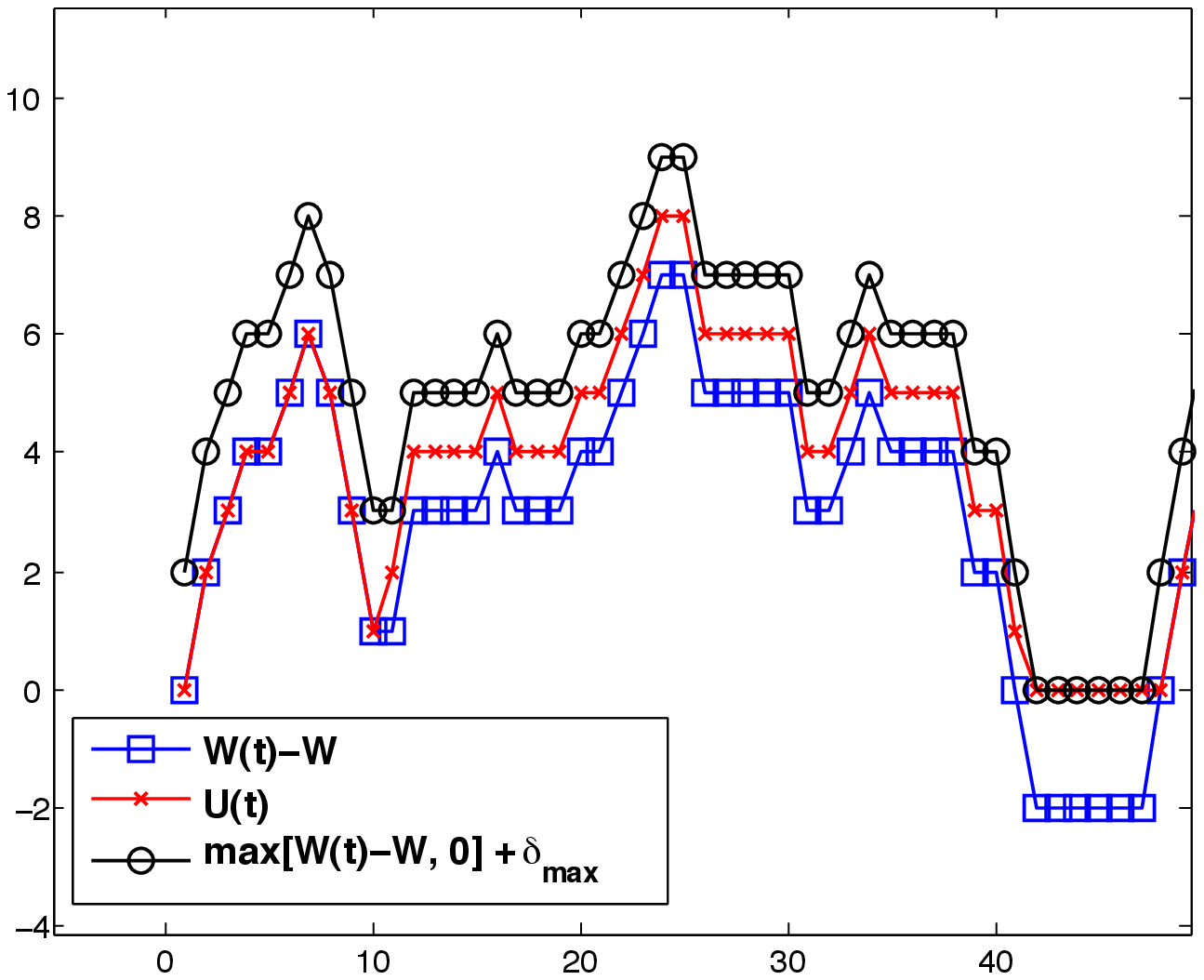}
\caption{Left: A sample backlog process; Right: An Example of $W(t)$ and $U(t)$.}\label{fig:samplequeue}
\end{figure}

To illustrate the idea of FQLA, we first look at an example. Figure \ref{fig:samplequeue} shows a $10^4$-slot sample backlog process under QLA.\footnote{This sample backlog process is one sample backlog process of queue $1$ of the system considered in Section \ref{section:numerical}, under QLA with $V=50$. } 
We see that after roughly 1500 slots, $U(t)$ always stays very close to $U^*_V$, which is a $\Theta(V)$ scalar in this case. To reduce delay,  we can first  find $\script{W}\in(0,U^*_V)$ such that: under QLA, there exists a time $t_0$ so that $U(t_0)\geq\script{W}$ and once $U(t)\geq\script{W}$, it remains so for all time (the solid line in Fig. \ref{fig:samplequeue} shows one for these $10^4$ slots). We then place $\script{W}$ fake bits (called \emph{place-holder bits} \cite{neelyrahul_asilomar}) in the queue at time $0$, i.e., initialize $U(0)=\script{W}$, and run QLA. It is easy to show that the utility performance of QLA will remain the same with this change, and the average backlog is now reduced by $\script{W}$. However, such a $\script{W}$ may require $\script{W}=U^*_V-\Theta(V)$, thus the average backlog may still be $\Theta(V)$. 

FQLA instead finds a $\script{W}$ such that in steady state, the backlog process under QLA \emph{rarely} goes below it, and places $\script{W}$ place-holder bits in the queue at time $0$.  FQLA  then uses an auxiliary process $W(t)$, called the \emph{virtual backlog process}, to keep track of the backlog process that should have been generated if QLA is used. Specifically, FQLA initializes $W(0)=\script{W}$. Then at every  slot, QLA is run using $W(t)$ as the queue size, and $W(t)$ is updated according to QLA. With $W(t)$ and $\script{W}$, FQLA works as follows: At time $t$, if $W(t)\geq \script{W}$, FQLA performs QLA's action (obtained based on $S(t)$ and $W(t)$); else if $W(t)<\script{W}$, FQLA carefully modifies QLA's action so as to maintain $U(t)\approx\max[W(t)-\script{W},0]$ for all $t$ (see Fig.\ref{fig:samplequeue} for an example). 
Similar as above, this roughly reduces the average backlog by $\script{W}$. The difference is that now we can show that $\script{W}=\max[U^*_V-\log^2(V),0]$ meets the requirement. Thus it is possible to bring the average backlog down to $O(\log^2(V))$. Also, since $W(t)$ can be viewed as a backlog process generated by QLA, it rarely goes below $\script{W}$ in steady state. Hence FQLA is almost always the same as QLA, thus is able to achieve an $O(1/V)$ close-to-optimal utility performance.

\subsection{The FQLA-Ideal Algorithm}
In this section, we present the FQLA-Ideal algorithm. We assume the value $\bv{U}^*_V=(U^*_{V1}, ..., U^*_{Vr})^T$ is known a-priori.

\emph{\underline{FQLA-Ideal:}} 
\begin{enumerate}
\item[(I)] \emph{\underline{Determining place-holder bits:}} 
For each $j$, define:  
\begin{eqnarray}
\script{W}_j=\max\big[U^*_{Vj}-\log^2(V),0\big], \label{eq:WLdef}
\end{eqnarray}
as the number of $\emph{place-holder bits}$ of queue $j$. 
\item[(II)] \emph{\underline{Place-holder-bit based action:}} Initialize
\[U_{j}(0)=0,\quad W_j(0)=\script{W}_j,  \quad\forall j.\]
For  $t\geq 1$, observe the network state $S(t)$, solve (\ref{eq:QLAeq}) with $\bv{W}(t)$ in place of $\bv{U}(t)$. Perform the chosen action with the following modification: Let $\bv{A}(t)$ and $\bv{\mu}(t)$ be the arrival and service rate vectors generated by the action. For each queue $j$, do (Idle fill whenever needed):
\begin{enumerate}
\item If $W_j(t)\geq \script{W}_j$: admit $A_j(t)$ arrivals, serve $\mu_j(t)$ data, i.e., update the backlog by: \[U_j(t+1)=\max\big[U_j(t)-\mu_j(t), 0\big]+A_j(t). \]

\item If $W_j(t)<\script{W}_j$: admit $\tilde{A}_j(t)=\max\big[A_j(t)-\script{W}_j+W_j(t), 0\big]$ arrivals, serve $\mu_j(t)$ data, i.e., update the backlog by:
\[U_j(t+1)=\max\big[U_j(t)-\mu_j(t),0\big]+\tilde{A}_j(t).\] 

\item Update $W_j(t)$ by:
\[W_j(t+1)=\max\big[W_j(t)-\mu_j(t), 0\big]+A_j(t). \]
\end{enumerate}
\end{enumerate}

From above we see that FQLA-Ideal is the same as QLA based on $\bv{W}(t)$ when $W_j(t)\geq \script{W}_j$ for all $j$. When $W_j(t)< \script{W}_j$ for some queue $j$, FQLA-Ideal admits roughly the \emph{excessive} packets after $W_j(t)$ is brought back to be above $\script{W}_j$ for the queue. Thus for problems where QLA admits an easy implementation, e.g., \cite{neelyenergy}, \cite{huangneelypricing}, it is also easy to implement FQLA. 
However,  we also notice two different features of FQLA: (1) By (\ref{eq:WLdef}), $\script{W}_j$ can be $0$. However, when $V$ is large, this happens only when $U^*_{0j}=U^*_{Vj}=0$ according to Lemma \ref{lemma:scalingeffectV}. In this case $\script{W}_j=U^*_{Vj}=0$, and queue $j$ indeed needs zero place-holder bits.  (2) Packets may be dropped in Step II-(b) upon their arrivals, or after they are admitted into the network in a multihop problem. Such packet dropping is natural in many flow control problems and does not change the nature of these problems. In other problems where such option is not available, the packet dropping option is introduced to achieve desired delay performance, and it can be shown that the fraction of packets dropped can be made arbitrarily small.  Note that packet dropping here is to compensate for the deviation from the desired Lagrange multiplier, thus is different from that in \cite{neelypktdrop}, where packet dropping is used for drift steering. 

\vspace{-.1in}
\subsection{Performance of FQLA-Ideal}
We look at the performance of FQLA-Ideal in this section. We first have the following lemma that shows the relationship between $\bv{U}(t)$ and $\bv{W}(t)$ under FQLA-Ideal. We will use it later to prove the delay bound of FQLA. Note that the lemma also holds for FQLA-General described later, as FQLA-Ideal/General differ only in the way of determining $\bv{\script{W}}=(\script{W}_1, ..., \script{W}_r)^T$.\begin{lemma}\label{lemma:UtWt}
Under FQLA-Ideal/General, we have $\forall\,\, j, t$:
\begin{equation}
\max\big[W_j(t)-\script{W}_j,0\big]\leq U_j(t)\leq\max\big[W_j(t)-\script{W}_j,0\big]+\delta_{max}\label{eq:UWrelation}
\end{equation}
where $\delta_{max}$ is defined in Section \ref{subsection:costtrafficservice} to be the upper bound of the number of arriving or departing packets of a queue. 
\end{lemma}
\begin{proof}
See Appendix C.
\end{proof}

The following theorem summarizes the main performance results of FQLA-Ideal. Recall that  for a given policy $\pi$, $f^{\pi}_{av}$ denotes its average cost defined in (\ref{eq:timeavcost}) and $f^{\pi}(t)$ denotes the cost induced by $\pi$ at time $t$. 

\begin{thm}\label{theorem:fqla_performance}
If the condition in Theorem \ref{theorem:prob_multicon_polyhedral} holds and a steady state distribution exists for the backlog process generated by QLA, then with a sufficiently large $V$, we have under FQLA-Ideal that,  
\begin{eqnarray}
\overline{U}&=&O(\log^2(V)),\label{theorem:fqla_backlog}\\
f^{FI}_{av}&=&f^*_{av}+O(1/V),\label{theorem:fqla_cost}\\
P_{drop} &=& O(1/V^{c_0\log(V)}),\label{theorem:fqla_pdrop}
\end{eqnarray}
where $c_0=\Theta(1)$, $\overline{U}$ is the time average backlog, $f_{av}^{FI}$ is the time average cost of FQLA-Ideal, $f^*_{av}$ is the optimal time average cost and $P_{drop}$ is the time average fraction of packets that are dropped in Step-II (b). 
\end{thm}
\begin{proof} Since a steady state distribution exists for the backlog process generated by QLA, we see that $\script{P}(D, m)$ in (\ref{eq:pm_def}) represents the steady state probability of the event that the backlog vector deviates from $\bv{U}^*_V$ by distance $D+m$. Now since $\bv{W}(t)$ can be viewed as a backlog process generated by QLA, with $\bv{W}(0)=\bv{\script{W}}$ instead of $\bv{0}$, we see from the proof of Theorem \ref{theorem:prob_multicon_polyhedral} that Theorem \ref{theorem:prob_multicon_polyhedral} and \ref{thm:prob_multi_con} hold for $\bv{W}(t)$, and by \cite{neelynowbook}, QLA based on $\bv{W}(t)$ achieves an average cost of $f^*_{av}+O(1/V)$. Hence by Theorem \ref{thm:prob_multi_con}, there exist constants $D_1, K_1, c_1^*=\Theta(1)$ so that:
$\script{P}^{(r)}(D_1, cK_1\log(V)) \leq \frac{c_1^*}{V^c}.$
By the definition of $\script{P}^{(r)}(D_1, cK_1\log(V))$, this implies that in steady state:
\[\prob\{W_j(t)>U^*_{Vj}+D_1+m\}\leq c_1^*e^{-\frac{m}{K_1}},\]
Now let: 
$Q_j(t)=\max[W_j(t)-U^*_{Vj}-D_1,0]$. We see that $\prob\{Q_j(t)>m\}\leq c_1^*e^{-\frac{m}{K_1}}$, $\forall m\geq0$.  We thus have $\overline{Q_j}= O(1)$, where $\overline{Q_j}$ is the time average value of $Q_j(t)$. Now it is easy to see by (\ref{eq:WLdef}) and (\ref{eq:UWrelation}) that $U_j(t)\leq Q_j(t)+\log^2(V)+D_1+\delta_{max}$ for all $t$. Thus (\ref{theorem:fqla_backlog}) follows since for a large $V$:
\begin{eqnarray*}
\overline{U_j}\leq \overline{Q_j}+\log^2(V)+D_1+\delta_{max} =\Theta(\log^2(V)).
\end{eqnarray*}

Now consider the average cost. To save space, we use FI for FQLA-Ideal. From above, we see that QLA based on $\bv{W}(t)$ achieves an average cost of $f^*_{av}+O(1/V)$. Thus it suffices to show that FQLA-Ideal performs almost the same as QLA based on $\bv{W}(t)$. First we have for all $t\geq1$ that:
\begin{eqnarray*}
\frac{1}{t}\sum_{\tau=0}^{t-1}f^{FI}(\tau)=\frac{1}{t}\sum_{\tau=0}^{t-1}f^{FI}(\tau)1_{E(\tau)}+\frac{1}{t}\sum_{\tau=0}^{t-1}f^{FI}(\tau)1_{E^c(\tau)}.
\end{eqnarray*}
Here $1_{E(\tau)}$ is the indicator function of the event $E(\tau)$, $E(\tau)$ is the event that FQLA-Ideal performs the same action as QLA at time $\tau$, and $1_{E^c(\tau)}=1-1_{E(\tau)}$. Taking expectation on both sides and using the fact that when FQLA-Ideal takes the same action as QLA, $f^{FI}(\tau)=f^{QLA}(\tau)$, we have:
\begin{eqnarray*}
\frac{1}{t}\sum_{\tau=0}^{t-1}\expect{f^{FI}(\tau)}&\leq&\frac{1}{t}\sum_{\tau=0}^{t-1}\expect{f^{QLA}(\tau)1_{E(\tau)}}\\
&&+\frac{1}{t}\sum_{\tau=0}^{t-1}\expect{\delta_{max}1_{E^c(\tau)}}.
\end{eqnarray*}
Taking the limit as $t$ goes to infinity on both sides and using $f^{QLA}(\tau)1_{E(\tau)}\leq f^{QLA}(\tau)$ , we get: 
\begin{eqnarray}
f^{FI}_{av}&\leq& f^{QLA}_{av} +\delta_{max}\lim_{t\rightarrow\infty}\frac{1}{t}\sum_{\tau=0}^{t-1}\expect{1_{E^c(\tau)}}\nonumber\\
&=& f^{QLA}_{av} +\delta_{max}\lim_{t\rightarrow\infty}\frac{1}{t}\sum_{\tau=0}^{t-1}\prob\{ E^c(\tau)\}.\label{eq:eq:fqla_av_cost}
\end{eqnarray}
However, $E^c(\tau)$ is included in the event that there exists a $j$ such that $W_j(\tau)<\script{W}_j$. Therefore by (\ref{eq:prob_pmr_special}) in Theorem \ref{thm:prob_multi_con}, for a large $V$ such that $\frac{1}{2}\log^2(V)\geq D_1$ and $\log(V)\geq8K_1$,
\begin{eqnarray}
\lim_{t\rightarrow\infty}\frac{1}{t}\sum_{\tau=0}^{t-1}\prob\{ E^c(\tau)\}&\leq& \script{P}^{(r)}(D_1, \log^2(V)-D_1) \nonumber\\
&=&O(c^*_1/V^{\frac{1}{2K_1}\log(V)})\nonumber\\
&=&O(1/V^4). \label{eq:fqla_erroreventbd}
\end{eqnarray}
Using this fact in (\ref{eq:eq:fqla_av_cost}), we obtain:
\begin{eqnarray*}
f^{FI}_{av}&=& f^{QLA}_{av} +O(\delta_{max}/V^4)= f^*_{av}+O(1/V), 
\end{eqnarray*}
where the last equality holds since $f^{QLA}_{av}=f^*_{av}+O(1/V)$. This proves (\ref{theorem:fqla_cost}). (\ref{theorem:fqla_pdrop}) follows since packets are dropped at time $\tau$ only if $E^c(\tau)$ happens, thus by (\ref{eq:fqla_erroreventbd}), the fraction of time when packet dropping happens is $O(1/V^{c_0\log(V)})$ with $c_0=\frac{1}{2K_1}=\Theta(1)$, and each time no more than $\sqrt{r}B$ packets can be dropped. 
\end{proof}



\subsection{The FQLA-General algorithm}
Now we describe the FQLA algorithm without any a-priori knowledge of $\bv{U}^*_V$, called FQLA-General. FQLA-General first runs the system for a long enough time $T$, such that the system enters its steady state. Then it chooses a sample of the queue vector value to estimate $\bv{U}^*_V$ and uses that to decide the number of place holder bits.  

\emph{\underline{FQLA-General:}} 
\begin{enumerate}
\item[(I)] \emph{\underline{Determining place-holder bits:}} 
\begin{enumerate}
\item Choose a large time $T$ (See Section \ref{section:fqla_practical} for the size of $T$)  and initialize $\bv{W}(0)=\bv{0}$. Run the QLA algorithm with parameter $V$, at every time slot $t$, update $\bv{W}(t)$ according to the QLA algorithm and obtain $\bv{W}(T)$. 
\item For each queue $j$, define:  
\begin{eqnarray}
\script{W}_j=\max\big[W_j(T)-\log^2(V),0\big],
\end{eqnarray}
as the number of $\emph{place-holder bits}$. 
\end{enumerate}
\item[(II)] \emph{\underline{Place-holder-bit based action:}} same as FQLA-Ideal.
\end{enumerate}
The performance of FQLA-General is summarized as follows:

\begin{thm}\label{theorem:fqla_performance2}
Assume the conditions in Theorem \ref{theorem:fqla_performance} hold and the system is in steady state at time $T$, then under FQLA-General with a sufficiently large $V$, with probability   $1-O(\frac{1}{V^4})$:
(a) $\overline{U}=O(\log^2(V))$, (b)
$f^{FG}_{av}=f^*_{av}+O(1/V)$, and (c) $P_{drop} = O(1/V^{c_0\log(V)})$, 
where $c_0=\Theta(1)$ and $f_{av}^{FG}$ is the time average cost of FQLA-General.
\end{thm}
\begin{proof} We will show that with probability of $1-O(\frac{1}{V^4})$, $\script{W}_j$ is close to $\max[U^*_{Vj}-\log^2(V), 0]$. The rest can then be proven similarly as in the proof of Theorem \ref{theorem:fqla_performance}. 

For each queue $j$, define:
\begin{eqnarray*}
v_j^+=U^*_{Vj}+\frac{1}{2}\log^2(V),\quad v^-_j=\max\big[U^*_{Vj}-\frac{1}{2}\log^2(V),0\big].
\end{eqnarray*}
Note that $v^-_j$ is defined with a $\max [\,]$ operator. This is due to the fact that $U^*_{Vj}$ can be zero. As in (\ref{eq:fqla_erroreventbd}), we see that by Theorem \ref{thm:prob_multi_con}, there exists $D_1=\Theta(1), K_1=\Theta(1)$ such that if $V$ is such that $\frac{1}{4}\log^2(V)\geq D_1$ and $\log(V)\geq16K_1$, then:
\begin{eqnarray*}
\prob\big\{\exists\, j,\,W_j(T)\notin[v_j^-, v_j^+]\big\}&\leq& \script{P}^{(r)}(D_1, \frac{1}{2}\log^2(V)-D_1)\\
&= &O(1/V^4)
\end{eqnarray*}
Thus we see that $\prob\big\{W_j(T)\in[v_j^-, v_j^+]\,\forall j\big\}=1- O(1/V^4)$, which implies: 
\begin{eqnarray*}
\prob\big\{\script{W}_j\in[\hat{v}_j^-, \hat{v}_j^+]\quad\forall j\big\}= 1- O(1/V^4).
\end{eqnarray*}
where $\hat{v}_j^+=\max\big[U^*_{Vj}-\frac{1}{2}\log^2(V), 0\big]$ and  $\hat{v}_j^-=\max\big[U^*_{Vj}-\frac{3}{2}\log^2(V),0\big]$. Hence for a large $V$, with probability  $1-O(\frac{1}{V^4})$: if $U^*_{Vj}>0$, we have $U^*_{Vj}-\frac{3}{2}\log^2(V)\leq \script{W}_j\leq U^*_{Vj}-\frac{1}{2}\log^2(V)$; else if $U^*_{Vj}=0$, we have $\script{W}_j=U^*_{Vj}$. The rest of the proof is similar as the proof of Theorem \ref{theorem:fqla_performance}.  
\end{proof}


\subsection{FQLA when $q_0()$ is locally smooth}\label{section:fqla_smooth}
Note that FQLA can also be implemented for problems with  $q_0(\bv{U})$ being locally smooth, with the only modification that $\script{W}_j=\max[U^*_{Vj}-\log^2(V)\sqrt{V},0]$. In this case, the following theorem can be obtained:
\begin{thm}
Assume the condition in Theorem \ref{theorem:qtwicediff} holds and a steady state distribution for a backlog process under QLA, then FQLA-Ideal achieves an $[O(1/V), O(\log^2(V)\sqrt{V})]$ performance-delay tradeoff, with $P_{drop}=O(1/V^{c_0\log(V)})$, where $c_0=\Theta(1)$; similarly, for appropriately chosen $T$, FQLA-General achieves the same performance with probability $1-O(1/V^4)$.
\end{thm}

\vspace{-.1in}
\subsection{Practical Issues}\label{section:fqla_practical}
From Lemma \ref{lemma:scalingeffectV} we see that the magnitude of $\bv{U}^*_V$ can be $\Theta(V)$. This means that $T$ in FQLA-General may need to be $\Omega(V)$, which is not very desirable when $V$ is large. We can instead use the following heuristic method to accelerate the process of determining $\bv{\script{W}}$: For every queue $j$, guess a very large $\script{W}_j$. Then start with this $\bv{\script{W}}$ and run the QLA algorithm for some $T_1$, say $\sqrt{V}$ slots. Observe the resulting backlog process. Modify the guess for each queue $j$ using a bisection algorithm until a proper $\bv{\script{W}}$ is found, i.e. when running QLA from that value, we observe fluctuations of $W_j(t)$ around $\script{W}_j$ instead of a nearly constant increase or decrease for all $j$. Then let $\script{W}_j=\max[\script{W}_j-\log^2(V),0]$ be the number of place-holder bits of queue $j$. To further reduce the error probability, one can repeat Step-I (a) multiple times and use the average value as $\bv{W}(T)$.

Note that even though results in Theorem \ref{theorem:fqla_performance} and \ref{theorem:fqla_performance2} assume a large $V$, in practice, the $V$ value may not have to be very large (See Section \ref{section:numerical} for an example).
 

\section{When there is a single queue}\label{section:secondspecial}
In this section, we look at the backlog process behavior under QLA under the special case when there is only one queue in the network. In this case, we have only a  single traffic constraint in the deterministic problem (\ref{eq:primal}):
\begin{eqnarray*}
\script{G}_1(\bv{x})=\sum_{s_i}p_{s_i}g_1(s_i, x^{(s_i)}) \leq\script{B}_1(\bv{x})=\sum_{s_i}p_{s_i}b_1(s_i, x^{(s_i)}),
\end{eqnarray*}
where $\bv{x}=(x^{(s_1)},...,x^{(s_M)})^T$. Thus $r=1$ and the Lagrange multiplier is a \emph{scalar}. This single queue setting is useful and can be used to model many network optimization problems, e.g., \cite{neelyenergy} and \cite{huangneelypricing}.  
Below, we first provide \emph{deterministic} upper and lower bounds for  $U(t)$. These bounds hold for \emph{arbitrary} network state distribution  and the way the state process evolves (possibly even non-ergodic). We then obtain a probabilistic bound of $U(t)$'s deviation from $U_V^*$ under general single queue network optimization problems. The probabilistic bound has the same form as those in Theorem \ref{theorem:prob_multicon_polyhedral} and \ref{theorem:qtwicediff}, but does not require any additional conditions such as (\ref{eq:multiple_dual_condition}) and (\ref{eq:qtwicediff_condition}).

\subsection{Deterministic Upper and Lower Bounds of $U(t)$}
Here we provide upper and lower bounds of $U(t)$ under QLA. First define the following problem for each network state $s_i$, for $i\in\{1,...,M\}$. 
\begin{eqnarray}
\max &&q_{s_i}(U)=\inf_{x^{(s_i)}\in\script{X}^{(s_i)}}\bigg\{Vf(s_i, x^{(s_i)})\label{eq:dual_isub}\\
&&\qquad\qquad\qquad\quad+U\big[g_{1}(s_i, x^{(s_i)})-b_1(s_i, x^{(s_i)})\big]\bigg\}\nonumber\\
s.t. && U\geq0.\nonumber
\end{eqnarray}
It is easy to see that $q_{s_i}(U)$ is the dual of (\ref{eq:primal}) when $s_i$ is the only network state. We now have the following theorem:

\begin{thm}\label{theorem:qla_u_ubbounds}
Assume (\ref{eq:dual_isub}) has a unique optimal solution $U^*_{s_i}\in[0,\infty]$ for all $s_i$. Consider the interval:  
\[\script{I}=\big[\min_{s_i}U^*_{s_i}-B, \max_{s_i}U^*_{s_i}+B\big],\]  
if under QLA, there exists $t_0\geq0$ such that $U(t_0)\in \script{I}$, then $U(t)\in \script{I}$ for all $t\geq t_0$.
\end{thm}

Note that here $[0,\infty]$ includes the value $\infty$. To prove Theorem \ref{theorem:qla_u_ubbounds}, we use the following lemma. 
\begin{lemma}\label{lemma:qla_tendtoopt}
If $U(t)\neq U^*_V$, then 
\begin{enumerate}
\item[(a)] Under QLA, \[\expect{(U(t)-U_V^*)\big[g_{1}(s_i, x^{(s_i)}_{U})-b_{1}(s_i, x^{(s_i)}_{U})\big]\left|\right.U(t)}<0.\]
\item[(b)] Under OSM, \[(U(t)-U_V^*)[\script{G}_1(x_U)-\script{B}_1(x_U)]<0.\]
\end{enumerate}
\end{lemma}
\begin{proof}
See Appendix D.
\end{proof}

Lemma \ref{lemma:qla_tendtoopt} shows that under QLA, if $U(t)<U_V^*$, then $\expect{g_1(s_i, x^{(s_i)}_{U})-b_1(s_i, x^{(s_i)}_{U})\left|\right. U(t)}>0$; else if  $U(t)>U_V^*$, we have $\expect{g_1(s_i, x^{(s_i)}_{U})-b_1(s_i, x^{(s_i)}_{U})\left|\right. U(t)}<0$. This shows that when $S(t)$ is i.i.d, the backlog value under QLA \emph{probabilistically} moves in the direction towards $U^*_V$. When there is a single network state, in which case (a) and (b) are equivalent, we see that $U(t)$ \emph{deterministically} moves in the direction towards $U^*_V$. 

\begin{proof} (Theorem \ref{theorem:qla_u_ubbounds})
First we see that, though it is possible for some $U^*_{s_i}$ to be infinity, it can be easily shown  that $\min_{s_i}U^*_{s_i}<\infty$. Thus $\script{I}$ is well defined. 

We now prove the lower bound. The upper bound can similarly be obtained. Without loss of generality, assume $U^*_{s_1}=\min_{s_i}U^*_{s_i}$ and $U^*_{s_M}=\max_{s_i}U^*_{s_i}$. Suppose at a time $t$ we have $U(t)\in \script{I}$:

(1) If  $U(t)\geq U^*_{s_1}$, we have $U(t+1)\geq U^*_{s_1}-B$, since $B$ is an upper bound of the magnitude change of $U(t)$. 

(2) Now if $U^*_{s_1}>U(t)\geq U^*_{s_1}-B$, we see that  $U(t)< U^*_{s_i}$ for all $i=1,...,M$. Also, when given $U(t)$ and $S(t)=s_i$, QLA's action is the same as OSM applied to (\ref{eq:dual_isub}). Thus by part (b) of Lemma \ref{lemma:qla_tendtoopt}, we see that $\script{G}_1(x_U)-\script{B}_1(x_U)=A(t)-\mu(t)>0$, hence by (\ref{eq:queuedynamic}) we have $U(t+1)>U(t)\geq U^*_{s_1}-B$. 
\end{proof}

Note that we did not use any assumption of the network state process in the above proof, hence   the result holds for \emph{arbitrary} network state distribution and the way $S(t)$ evolves. 

\subsection{Probabilistic bound of $U(t)$'s deviation from $U^*_V$}
In this section we provide a probabilistic bound of $U(t)$'s deviation from $U_V^*$. The bound has a similar form as those in Theorem \ref{theorem:prob_multicon_polyhedral} and \ref{theorem:qtwicediff}, but only applies to general single queue optimization problems. However, the bound here does not require additional conditions such as (\ref{eq:multiple_dual_condition}) and (\ref{eq:qtwicediff_condition}). Hence it is more general than the previous results when restricted to single queue optimization problems. 
Recall that $\script{P}(D,m)$ is defined in (\ref{eq:pm_def}) as:
\begin{eqnarray*}
\script{P}(D, m)\triangleq\lim\sup_{t\rightarrow\infty}\frac{1}{t}\sum_{\tau=0}^{t-1}\prob\{|U(\tau)-U_V^*|>D+m\}.
\end{eqnarray*}

\begin{thm}\label{theorem:prob_singlecon}
Under QLA, there exist constants $d,a^*, \rho^*>0$, possibly dependent on $V$, such that:
\begin{eqnarray}
\script{P}(d,m)\leq a^* e^{-\rho^*m}.\label{eq:pm_ineq_r1}
\end{eqnarray}
\end{thm}

Theorem \ref{theorem:prob_singlecon} shows that the probability to deviate further from $U^*_V$ will eventually be exponential. To prove Theorem \ref{theorem:prob_singlecon}, we need the following lemmas: 
\begin{lemma}\label{lemma:dfinite}
$q(\bv{U}^*_V)>-\infty$. 
\end{lemma}
\begin{lemma}\label{lemma:moretendtoopt}
Under QLA, if (a) $0\leq U_1<U_2<U_V^*$ or (b) $0\leq U_V^*<U_1<U_2$, then:
\begin{eqnarray*}
\expect{g_{1}(s_i, x^{(s_i)}_{U_1})-b_{1}(s_i, x^{(s_i)}_{U_1})\left|\right.U_1(t)}\qquad\qquad\qquad\qquad\\
\geq\expect{ g_{1}(s_i, x^{(s_i)}_{U_2})-b_{1}(s_i, x^{(s_i)}_{U_2})\left|\right.U_2(t)}.
\end{eqnarray*}
In case (a), both quantities are positive; while in case (b), both quantities are negative. 
\end{lemma}
\begin{lemma}\label{lemma:subgradient_qla}
Under QLA,
\begin{eqnarray}
(U_{V}^*-U(t))\expect{\big[g_1(s_i, x^{(s_i)})-b_1(s_i, x^{(s_i)})\big]\left|\right. U(t)}.\label{eq:qla_osm_relation}\\
\geq q(U^*_V)-q(U(t))\nonumber
\end{eqnarray}
\end{lemma}

Lemma \ref{lemma:dfinite} follows  easily from the $\epsilon$-slackness assumption in Section \ref{subsection:costtrafficservice}. 
Lemma \ref{lemma:moretendtoopt} can be viewed as saying that when $U(t)$ deviates more from $U^*_V$, the chosen action generates a larger drift  towards $U^*_V$. Lemma \ref{lemma:subgradient_qla} can be viewed as the subgradient property under QLA. Lemma \ref{lemma:moretendtoopt} and \ref{lemma:subgradient_qla} are proven in Appendix E. We now take the following approach to prove Theorem \ref{theorem:prob_singlecon}. We first use Lemma \ref{lemma:dfinite} and \ref{lemma:subgradient_qla} to find a single $U(t)$ value, whose drift value is large enough for analysis, and then conclude by Lemma  \ref{lemma:moretendtoopt} that any other $U(t)$ that is further away from $U^*_V$ generates a larger drift. Then we carry out the same drift analysis as in the proof of Theorem \ref{theorem:prob_multicon_polyhedral} to obtain the probability bound. 

\begin{proof} (Theorem \ref{theorem:prob_singlecon})
Since $r=1$, we have the dual function being:
\begin{eqnarray*}
q(U)=\inf_{x^{(s_i)}\in \script{X}^{(s_i)}}\bigg\{ V\sum_{s_i}p_{s_i}f(s_i, x^{(s_i)})\qquad\qquad\qquad\\
+U\big[ \sum_{s_i}p_{s_i}g_1(s_i, x^{(s_i)})- \sum_{s_i}p_{s_i}b_1(s_i, x^{(s_i)})\big]\bigg\}.
\end{eqnarray*}
Now by the $\epsilon$-slackness assumption in Section \ref{subsection:costtrafficservice} and  the fact that the cost functions are bounded by $\delta_{max}$, it can easily be shown that: 
\begin{eqnarray*}
q(U)\leq V\delta_{max} - \epsilon U\quad\forall\,\,U\geq0. 
\end{eqnarray*}
Hence if $q(U)\geq q(U^*_V)-\epsilon_0$ for some $\epsilon_0\geq0$, then we have:
\begin{eqnarray*}
\epsilon_0\geq q(U^*_V)-q(U)\geq q(U^*_V)+\epsilon U-V\delta_{max},
\end{eqnarray*}
which by Lemma \ref{lemma:dfinite} implies:
\begin{eqnarray}
U\leq \frac{\epsilon_0+V\delta_{max}-q(U^*_V)}{\epsilon}<\infty.\label{eq:aUsizebound}
\end{eqnarray}

Now fix an $\epsilon_0>0$, define the set $S_{\epsilon_0}=\{U\geq 0\, |\, q(U) \geq q(U_V^*)-\epsilon_0\}$. Define: 
\begin{eqnarray}
d(V,\epsilon_0)=\sup_{U\in S_{\epsilon_0}} |U-U_V^*|. \label{eq:dv_def}
\end{eqnarray}
By (\ref{eq:aUsizebound}) we see that $d(V,\epsilon_0)\in(0, \infty)$. Also whenever $|U(t)-U^*_V|>d(V,\epsilon_0)$, we have:
\begin{eqnarray}
q(U_V^*)-q(U(t))\geq\epsilon_0.
\end{eqnarray}
Thus by Lemma \ref{lemma:subgradient_qla}, we see that when $|U(t)-U^*_V|>d(V,\epsilon_0)$,
\begin{eqnarray*}
(U^*_V-U(t))\expect{(g_1(s_i, x^{(s_i)}_{U})-b_1(s_i, x^{(s_i)}_{U}) )\left|\right.U(t)}
\geq\epsilon_0.
\end{eqnarray*}

Now consider $U^*_V>d(V,\epsilon_0)+\epsilon_1$ for some small $\epsilon_1>0$. Define $U_l\triangleq U^*_V-d(V,\epsilon_0)-\epsilon_1$. From above and Lemma \ref{lemma:subgradient_qla} we see that if $U(t)=U_l$, then:
\begin{eqnarray}
\expect{(g_1(s_i, x^{(s_i)}_{U})-b_1(s_i, x^{(s_i)}_{U}) )\left|\right.U(t)}\geq\frac{\epsilon_0}{d(V,\epsilon_0)+\epsilon_1}.\label{eq:driftlowerbound_singlecon}
\end{eqnarray}
Denote $\eta_d = \frac{\epsilon_0}{d(V,\epsilon_0)+\epsilon_1}$. It is easy to see by (\ref{eq:Bdef}) that $\eta_d\leq B$. Using Lemma \ref{lemma:moretendtoopt}, we see that (\ref{eq:driftlowerbound_singlecon}) holds for all $U(t)\leq U_l= U^*_V-d(V,\epsilon_0)-\epsilon_1$. A similar argument will show that whenever $U(t)\geq U_u\triangleq U^*_V+d(V,\epsilon_0)+\epsilon_1$, 
\begin{eqnarray}
\expect{(g_1(s_i, x^{(s_i)}_{U})-b_1(s_i, x^{(s_i)}_{U}) )\left|\right.U(t)}\leq-\eta_d.\label{eq:driftlowerbound2_singlecon}
\end{eqnarray}
Now let $d=d(V,\epsilon_0)+\epsilon_1$ and define:
\begin{eqnarray}
Y(t)=\max\{|U(t)-U^*_V|-d, 0\},\label{eq:virtualprocess}
\end{eqnarray}
then whenever $Y(t)\geq B$, we have \[\expect{Y(t+1)-Y(t)\left|\right.U(t)}\leq-\eta_d.\] 
Also $|Y(t+1)-Y(t)|\leq B$ for all $t$. We can now carry out a similar argument as in the proof of Theorem \ref{theorem:prob_multicon_polyhedral} and obtain: 
\begin{eqnarray}
\lim\sup_{t\rightarrow\infty}\frac{1}{t}\sum_{\tau=1}^t\frac{w\eta_d}{2}e^{wm}\prob\{Y(\tau)>m\} &\leq&e^{2wB},
\end{eqnarray}
where $w=\frac{\eta_d}{B^2+B\eta_d/3}$. Thus we have:
\begin{eqnarray}
\script{P}(d, m) \leq \frac{2(B^2+B\eta_d/3)e^{\frac{2\eta_d}{B+\eta_d/3}}}{\eta_d^2} e^{-\frac{\eta_d m}{B^2+B\eta_d/3}}. \label{eq:singlecon_devbound_foo}
\end{eqnarray}
Therefore (\ref{eq:pm_ineq_r1}) holds with: 
\begin{eqnarray}
a^*=\frac{2(B^2+B\eta_d/3)e^{\frac{2\eta_d}{B+\eta_d/3}}}{\eta_d^2},\quad  \rho^*=\frac{\eta_d}{B^2+B\eta_d/3}.\label{eq:pm_ineq_const_r1}
\end{eqnarray}

Now if $U^*_V-d(V,\epsilon_0)-\epsilon_1<0$, then we have $U^*_V-U(t)\leq d$ whenever $U(t)\leq U^*_V$. Thus the $\{Y(\tau)>m\}$ is simply the event that $U(t)>U^*_V+d+m$. It is easy to see from above that (\ref{eq:singlecon_devbound_foo}) also holds in this case. 
\end{proof}

To see how Theorem \ref{theorem:prob_singlecon} is related to Theorem \ref{theorem:prob_multicon_polyhedral} and \ref{theorem:qtwicediff},    first consider (\ref{eq:multiple_dual_condition}) holds for all $U\geq0$. In this case, for a fixed $\epsilon_0=\Theta(1)$, we  have for all $U\in S_{\epsilon_0}$ that: 
\begin{eqnarray*}
\epsilon_0\geq q(U^*_V)-q(U)\geq L|U-U^*_V|.
\end{eqnarray*}
Thus $d(V, \epsilon_0)=\Theta(1)$, which then implies $\eta_d, \rho^*$ and $a^*$ are all $\Theta(1)$. Thus by (\ref{eq:singlecon_devbound_foo}) we see that $U(t)$ will mostly be within $O(\log(V))$ distance from $U^*_V$, as stated in Theorem \ref{theorem:prob_multicon_polyhedral}. Now if (\ref{eq:qtwicediff_condition}) holds for all $U\geq0$, then we see from (\ref{eq:qcon_twicediff_generalV}) that: 
\begin{eqnarray*}
\epsilon_0\geq q(U^*_V)-q(U)\geq \frac{L}{V}|U-U^*_V|^2 \quad\forall\,\, U\in S_{\epsilon_0}.
\end{eqnarray*}
This implies $d(V,\epsilon_0)=O(\sqrt{V})$ and $\eta_d=\Omega(1/\sqrt{V})$. Thus $\rho^*=\Omega(1/\sqrt{V})$ and $a^*=O(V)$ and again $U(t)$ is mostly within $O(\sqrt{V}\log(V))$ distance from $U^*_V$, as shown in Theorem \ref{theorem:qtwicediff}. 

\vspace{-.1in}
\section{Simulation}\label{section:numerical}
\begin{figure}[cht]
\centering
\includegraphics[height=0.6in, width=3.65in]{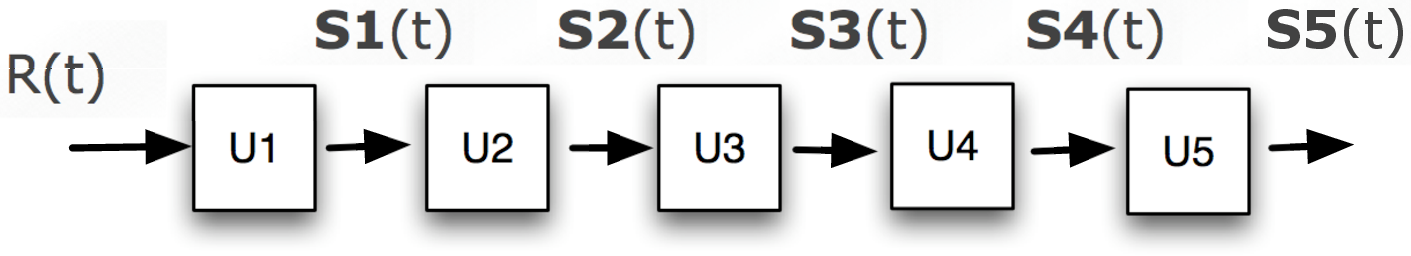}
\caption{A five queue system}\label{fig:EECA_5q}
\end{figure}

In this section we provide simulation results for the FQLA algorithms. For simplicity, we only consider the case where $q_0(\bv{U})$ is locally polyhedral. We consider a five queue system that extends the example in Section \ref{subsection:toyexample}. In this case  $r=5$. The system is shown in Fig. \ref{fig:EECA_5q}. The goal is to perform power allocation at each node so as to support the arrival with minimum energy expenditure. 

In this example, the random network state $S(t)$ is the vector containing the random arrivals $R(t)$ and the channel states $S_i(t)$, $i=1,...,5$. Similar as in Section \ref{subsection:toyexample}, we have:
\begin{eqnarray*}
\bv{A}(t)&=&(R(t), \mu_1(t), \mu_2(t), \mu_3(t), \mu_4(t))^T,\\
\bv{\mu}(t)&=&( \mu_1(t), \mu_2(t), \mu_3(t), \mu_4(t), \mu_5(t))^T,
\end{eqnarray*}
i.e., $A_1(t)=R(t)$, $A_i(t)=\mu_{i-1}(t)$ for $i\geq2$, where $\mu_i(t)$ is the service rate obtained by queue $i$ at time $t$.
$R(t)$ is $0$ or $2$ with probabilities $\frac{3}{8}$ and $\frac{5}{8}$, respectively. $S_i(t)$ can be ``Good'' or ``Bad'' with equal probabilities for $1\leq i\leq 5$. When the channel is good, one unit of power can serve two packets; otherwise one unit of power can serve only one packet. We assume all channels can be activated at the same time without affecting others. It can be verified that $\bv{U}^*_V=(5V,4V,3V,2V,V)^T$ is unique. In this example, the backlog vector process evolves as a Markov chain with countably many states. Thus one can show that there exists a stationary distribution for the backlog vector under QLA. 
  
\begin{figure}[cht]
\centering
\includegraphics[height=1.42in, width=2.2in]{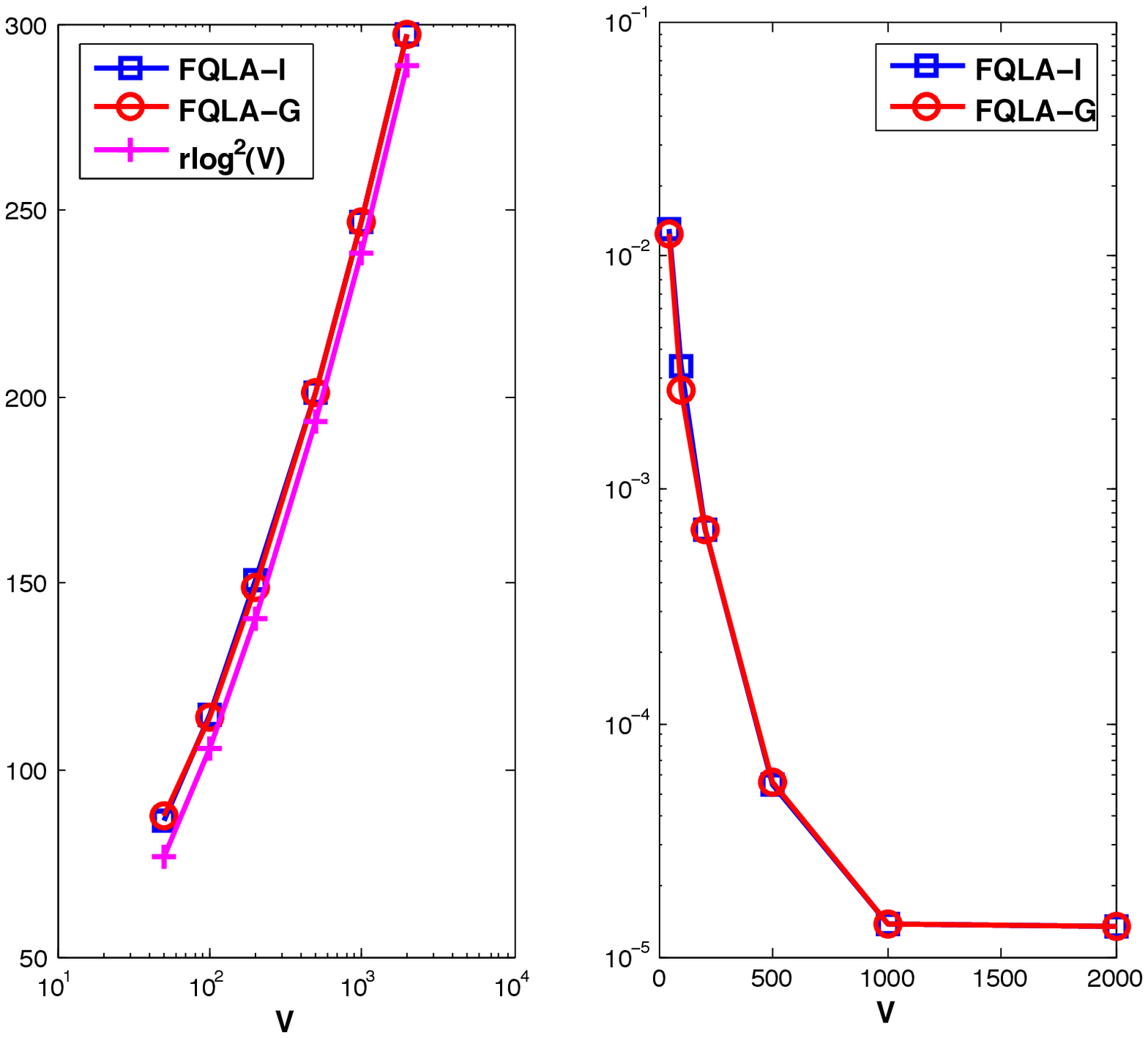}
\includegraphics[height=1.5in, width=1.22in]{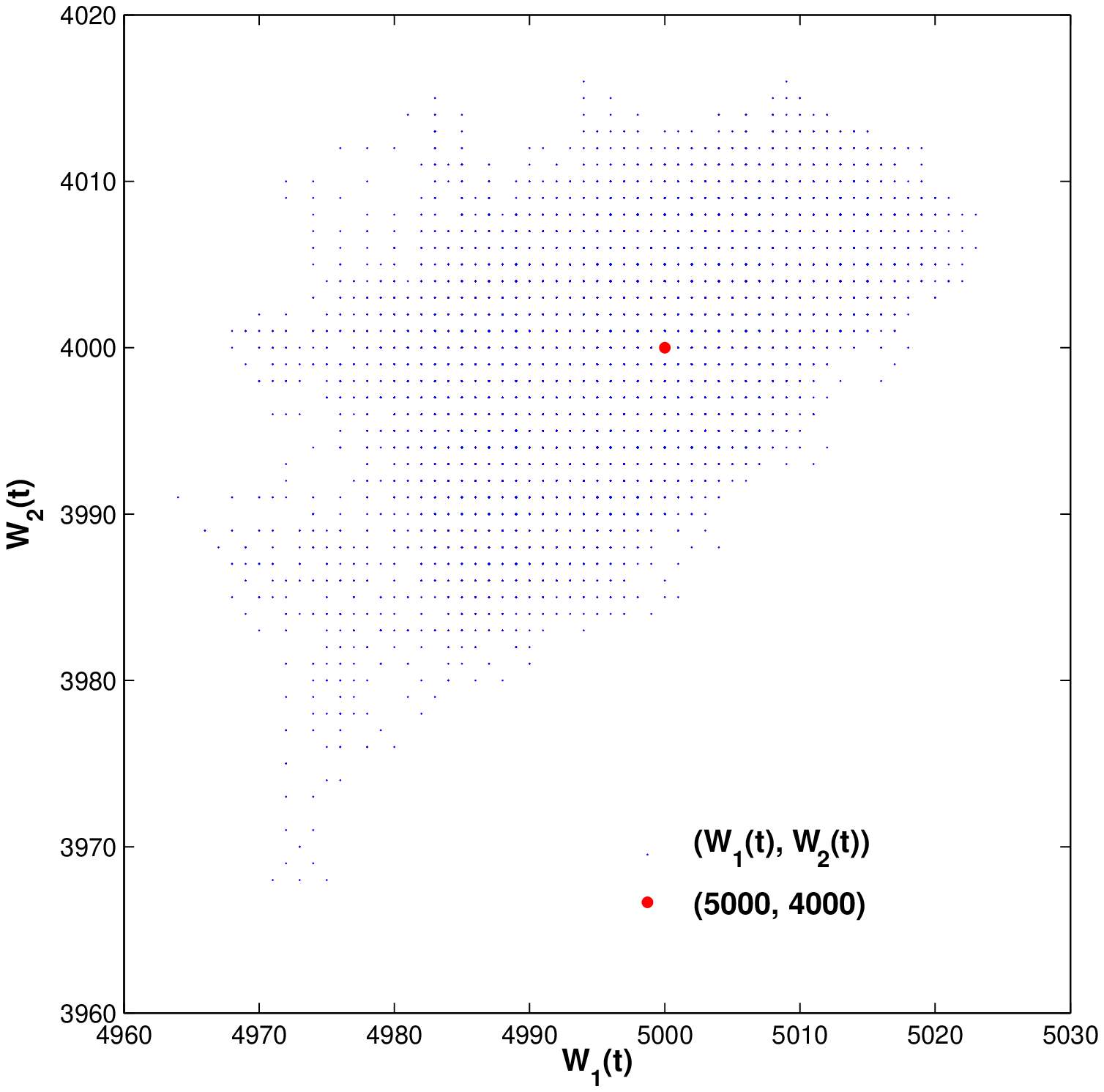}
\caption{FQLA-Ideal performance: Left - Average queue size; Middle - Percentage of packets dropped; Right - Sample $(W_1(t), W_2(t))$ process for $t\in[10000,110000]$ and $V=1000$ under FQLA-Ideal.}
\label{fig:FQLA result show}
\end{figure}

 We simulate FQLA-Ideal and FQLA-General with $V=50, 100, 200, 500, 1000$ and $2000$. We run each case for $5\times10^6$ slots under both algorithms. For FQLA-General, we use $T=50V$ in Step-I and  repeat Step-I $100$ times  and use their average as $\bv{W}(T)$. 
It is easy to see from the left plot in Fig. \ref{fig:FQLA result show} that the average queue sizes under both FQLAs are always close to the value $5\log^2(V)$ ($r=5$). From the middle plot we also see that the percentage of packets dropped decreases rapidly and gets below $10^{-4}$ when $V\geq500$ under both FQLAs. These plots show that in practice, $V$ may not have to be very large for Theorem \ref{theorem:fqla_performance} and \ref{theorem:fqla_performance2} to hold. The right plot shows a sample $(W_1(t), W_2(t))$ process for a $10^5$-slot interval under FQLA-Ideal with $V=1000$, considering only the first two queues of Fig. \ref{fig:EECA_5q} for this example. We see that during this interval, $(W_1(t),W_2(t))$ always remains close to $(U^*_{V1}, U^*_{V2})=(5V,4V)$, and $W_1(t)\geq\script{W}_1=4952$, $W_2(t)\geq\script{W}_2=3952$. 
For all $V$ values, the average power expenditure is very close to $3.75$, which is the optimal energy expenditure, and the average of $\sum W_j(t)$ is very close to $15V$ (plots  omitted for brevity). 
 
\vspace{-.25in}
\section{Lagrange Multiplier: ``shadow price'' and ``network gravity''}\label{section:LMfunctionality}It is well known that Lagrange Multipliers can play the role of ``shadow prices'' to regulate flows in many flow-based problems with different objectives, e.g., \cite{kellyelastic97}. This important feature has enabled the development of many distributed algorithms in resource allocation problems, e.g., \cite{curescupriceutility05}. However, a problem of this type typically requires data transmissions to be represented as flows. Thus in a network that is discrete in nature, e.g., time slotted or packetized transmission, a rate allocation solution obtained by solving such a flow-based problem does not immediately specify a scheduling policy. 

Recently, several Lyapunov algorithms have been proposed to solve utility optimization problems under discrete network settings. In these algorithms, backlog vectors act as the ``gravity'' of the network and allow optimal scheduling to be built upon them. It is also revealed in \cite{neelypowerjsac} that QLA is closely related to the dual subgradient method and backlogs play the same role as Lagrange multipliers in a time invariant network. Now we see by Theorem \ref{theorem:prob_multicon_polyhedral} and \ref{theorem:qtwicediff}  that the backlogs indeed play the same role as Lagrange multipliers even under a more general stochastic network. 

In fact, the backlog process under QLA can be closely related to a sequence of updated Lagrange multipliers under a subgradient method. Consider the following important variant of OSM, called  the randomized incremental subgradient method (RISM) \cite{bertsekasoptbook}, which makes use of the separable nature of (\ref{eq:dual_separable}) and solves the dual problem (\ref{eq:dualproblem}) as follows:

\underline{\emph{RISM:}} Initialize $\bv{U}(0)$; at iteration $t$, observe $\bv{U}(t)$, choose a random state $S(t)\in\script{S}$ according to some probability law.
(1) If $S(t)=s_i$, find $x_{\bv{U}}^{(s_i)}\in\script{X}^{(s_i)}$ that solves the following:
\begin{eqnarray}
\hspace{-.1in}\min && Vf(s_i, x)+\sum_jU_j(t)\big[g_j(s_i, x)-b_j(s_i, x)\big]\nonumber \\
\hspace{-.1in}s.t. && x\in \script{X}^{(s_i)}.\label{eq:dual_ith_sub}
\end{eqnarray}
(2) Using the $x_{\bv{U}}^{(s_i)}$ found, update $\bv{U}(t)$ according to:  \footnote{Note that this update rule is different from RISM's usual rule, i.e., $U_j(t+1)=\max\big[U_j(t)-\alpha^t b_j(s_i, x)+\alpha^t g_j(s_i, x), 0\big]$, but it almost does not affect the performance of RISM.}
\begin{eqnarray*}
\hspace{-.05in}U_j(t+1)=\max\bigg[U_j(t)-\alpha^t b_j(s_i, x_{\bv{U}}^{(s_i)}),0\bigg]+\alpha^t g_j(s_i, x_{\bv{U}}^{(s_i)}).
\end{eqnarray*}

As an example, $S(t)$ can be chosen by independently choosing $S(t)=s_i$ with probability $p_{s_i}$ every time slot. In this case $S(t)$ will be i.i.d..
Note that in the stochastic problem, a network state $s_i$ is chosen randomly by nature as the physical system state at time $t$; while here a state is chosen artificially by RISM according some probability law. Now we see from (\ref{eq:QLAeq}) and (\ref{eq:dual_ith_sub}) that: given the same $\bv{U}(t)$ and $s_i$, \emph{QLA and RISM choose an action in the same way}.  If also $\alpha^t=1$ for all $t$, and that $S(t)$ under RISM evolves according to the same probability law as $S(t)$ of the physical system, we see that \emph{applying QLA to the network is indeed equivalent to applying RISM to the dual problem of (\ref{eq:primal}), with the network state being chosen by nature, and the network backlog being the Lagrange multiplier.} Therefore,  Lagrange Multipliers under such stochastic discrete network settings act as the ``network gravity,'' thus allow scheduling to be done optimally and adaptively based on them. This ``network gravity'' functionality of Lagrange Multipliers in discrete network problems can thus be viewed as the counterpart of their ``shadow price'' functionality in the flow-based problems. Further more, the ``network gravity'' property of Lagrange Multipliers enables the use of place holder bits to reduce network delay in network utility optimization problems. This is a unique feature not possessed by its ``price'' counterpart.

\section*{Appendix A- Proof of Lemma \ref{lemma:expecteddistance}}
Here we prove Lemma \ref{lemma:expecteddistance}. First we prove the following useful lemma. 
\begin{lemma}\label{lemma:queuedynamic_diff}
Under queueing dynamic (\ref{eq:queuedynamic}), we have: 
\begin{eqnarray*}
\|\bv{U}(t+1)-\bv{U}^*_V\|^2&\leq& \|\bv{U}(t)-\bv{U}^*_V\|^2+2B^2\\
&&- 2\big(\bv{U}^*_V-\bv{U}(t)\big)^T(\bv{A}(t)-\bv{\mu}(t)).
\end{eqnarray*}
\end{lemma}
\begin{proof} (Lemma \ref{lemma:queuedynamic_diff})
From (\ref{eq:queuedynamic}), we see that $\bv{U}(t+1)$ is obtained by first projecting $\bv{U}(t)-\bv{\mu}(t)$ onto $\mathbb{R}^r_+$ and then adding $\bv{A}(t)$. Thus we have (we use $[\bv{x}]^+$ to denote the projection of $\bv{x}$ onto $\mathbb{R}^r_+$): 
\begin{eqnarray}
&&\|\bv{U}(t+1)-\bv{U}^*_V\|^2\nonumber\\
&=&\|[\bv{U}(t)-\bv{\mu}(t)]^++\bv{A}(t)-\bv{U}^*_V\|^2\nonumber\\
&=&\big([\bv{U}(t)-\bv{\mu}(t)]^++\bv{A}(t)-\bv{U}^*_V\big)^T\nonumber\\
&&\qquad\qquad\qquad\qquad\big([\bv{U}(t)-\bv{\mu}(t)]^++\bv{A}(t)-\bv{U}^*_V\big)\nonumber\\
&=&\big([\bv{U}(t)-\bv{\mu}(t)]^+-\bv{U}^*_V\big)^T\big([\bv{U}(t)-\bv{\mu}(t)]^+-\bv{U}^*_V\big)\nonumber\\
&& + 2\big([\bv{U}(t)-\bv{\mu}(t)]^+-\bv{U}^*_V\big)^T\bv{A}(t)+\|\bv{A}(t)\|^2.\label{eq:lemmadiffq_eq1}
\end{eqnarray}
Now by the non expansive property of projection  \cite{bertsekasoptbook}, we have:
\begin{eqnarray*}
&&\big([\bv{U}(t)-\bv{\mu}(t)]^+-\bv{U}^*_V\big)^T\big([\bv{U}(t)-\bv{\mu}(t)]^+-\bv{U}^*_V\big)\\
&\leq& \big(\bv{U}(t)-\bv{\mu}(t)-\bv{U}^*_V\big)^T\big(\bv{U}(t)-\bv{\mu}(t)-\bv{U}^*_V\big)\\
&=&\|\bv{U}(t)-\bv{U}^*_V\|^2+\|\bv{\mu}(t)\|^2-2(\bv{U}(t)-\bv{U}^*_V)^T\bv{\mu}(t).
\end{eqnarray*}
Plug this into (\ref{eq:lemmadiffq_eq1}), we have:
\begin{eqnarray}
\hspace{-.3in}&&\|\bv{U}(t+1)-\bv{U}^*_V\|^2\label{eq:lemmadiffq_eq2}\\
\hspace{-.3in}&\leq& \|\bv{U}(t)-\bv{U}^*_V\|^2+\|\bv{\mu}(t)\|^2-2(\bv{U}(t)-\bv{U}^*_V)^T\bv{\mu}(t)\nonumber\\
\hspace{-.3in}&& + \|\bv{A}(t)\|^2+ 2\big([\bv{U}(t)-\bv{\mu}(t)]^+-\bv{U}^*_V\big)^T\bv{A}(t).\nonumber
\end{eqnarray}
Now since $\bv{U}(t), \bv{\mu}(t), \bv{A}(t)\succeq\bv{0}$,  it is easy to see that:
\begin{eqnarray}
\big([\bv{U}(t)-\bv{\mu}(t)]^+\big)^T\bv{A}(t)\leq \bv{U}(t)^T\bv{A}(t).\label{eq:lemmadiffq_eq3}
\end{eqnarray}
By (\ref{eq:lemmadiffq_eq2}) and (\ref{eq:lemmadiffq_eq3}) we have:
\begin{eqnarray*}
\hspace{-.2in}&&\|\bv{U}(t+1)-\bv{U}^*_V\|^2\\
\hspace{-.2in}&\leq& \|\bv{U}(t)-\bv{U}^*_V\|^2+\|\bv{\mu}(t)\|^2-2(\bv{U}(t)-\bv{U}^*_V)^T\bv{\mu}(t)\\
\hspace{-.2in}&& + \|\bv{A}(t)\|^2+ 2\big(\bv{U}(t)-\bv{U}^*_V\big)^T\bv{A}(t)\\
\hspace{-.2in}&\leq& \|\bv{U}(t)-\bv{U}^*_V\|^2+2B^2- 2\big(\bv{U}^*_V-\bv{U}(t)\big)^T(\bv{A}(t)-\bv{\mu}(t)), 
\end{eqnarray*}
where the last inequality follows since $\|\bv{A}(t)\|^2\leq B^2$ and $\|\bv{\mu}(t)\|^2\leq B^2$.
\end{proof}

We now prove Lemma \ref{lemma:expecteddistance}.

\begin{proof} (Lemma \ref{lemma:expecteddistance})
By Lemma \ref{lemma:queuedynamic_diff} we see that when $S(t)=s_i$, we have the following for any network state $s_i$ with a given $\bv{U}(t)$ (here we add superscripts to $\bv{U}(t+1)$, $\bv{A}(t)$ and $\bv{\mu}(t)$ to indicate their dependence on $s_i$):
\begin{eqnarray}
\|\bv{U}^{(s_i)}(t+1)-\bv{U}_V^*\|^2\leq \|\bv{U}(t)-\bv{U}_V^*\|^2+2B^2\qquad\quad\label{eq:qla_iteration_isub_si}\\
\qquad\qquad-2(\bv{U}_V^*-\bv{U}(t))^{T}(\bv{A}^{(s_i)}(t)-\bv{\mu}^{(s_i)}(t)).\nonumber
\end{eqnarray}
By definition, $A^{(s_i)}_j(t)=g_j(s_i, x^{(s_i)}_{\bv{U}})$, and $\mu^{(s_i)}_j(t)=b_j(s_i, x^{(s_i)}_{\bv{U}})$, with $x^{(s_i)}_{\bv{U}}$ being the solution of (\ref{eq:QLAeq}) for the given $\bv{U}(t)$. Now  consider the deterministic problem (\ref{eq:primal}) with only a single network state $s_i$, then the corresponding dual function (\ref{eq:dual}) becomes:
\begin{eqnarray}
q_{s_i}(\bv{U}(t))=\inf_{x^{(s_i)}\in\script{X}^{(s_i)}}\bigg\{Vf(s_i, x^{(s_i)})\label{eq:dual_singlestate}\qquad\qquad\qquad\qquad\\
+\sum_jU_j(t)\big[g_j(s_i, x^{(s_i)})-b_j(s_i, x^{(s_i)})\big]\bigg\}.\nonumber
\end{eqnarray}
Therefore by (\ref{eq:subgradient_def}) we see that $(\bv{A}^{(s_i)}(t)-\bv{\mu}^{(s_i)}(t))$ is a subgradient of $q_{s_i}(\bv{U})$ at $\bv{U}(t)$. Thus by (\ref{eq:subgradientpro}) we have:
\begin{eqnarray}
(\bv{U}_V^*-\bv{U}(t))^{T}(\bv{A}^{(s_i)}(t)-\bv{\mu}^{(s_i)}(t))\label{eq:qla_osm_relation_pre}\qquad\qquad\qquad\\
\geq q_{s_i}(\bv{U^*_V})-q_{s_i}(\bv{U}(t)).\nonumber
\end{eqnarray}
Plug (\ref{eq:qla_osm_relation_pre}) into (\ref{eq:qla_iteration_isub_si}), we get:
\begin{eqnarray}
\|\bv{U}^{(s_i)}(t+1)-\bv{U}_V^*\|^2&\leq& \|\bv{U}(t)-\bv{U}_V^*\|^2+2B^2\label{eq:qla_iteration_dist}\qquad\qquad\\
\qquad&-&2\big(q_{s_i}(\bv{U^*_V})-q_{s_i}(\bv{U}(t))\big).\nonumber
\end{eqnarray}
More generally, we have:
\begin{eqnarray}
\|\bv{U}(t+1)-\bv{U}^*_V\|^2 \leq \|\bv{U}(t)-\bv{U}^*_V\|^2 + 2B^2\label{eq:qla_iteration_general}\qquad\qquad \\
-2\big(q_{S(t)}(\bv{U}^*_V)-q_{S(t)}(\bv{U}(t))\big). \nonumber
\end{eqnarray}
Now fix $\nu>0$, summing up (\ref{eq:qla_iteration_general}) from time $t$ to $t+T_{\nu}-1$, we obtain:
\begin{eqnarray}
\|\bv{U}(t+T_{\nu})-\bv{U}^*_V\|^2\leq \|\bv{U}(t)-\bv{U}^*_V\|^2 + 2T_{\nu}B^2\label{eq:Tslotdrift_pre}\qquad\quad\\  
-2\sum_{\tau=0}^{T_{\nu}-1}\big[q_{S(t+\tau)}(\bv{U}^*_V)-q_{S(t+\tau)}(\bv{U}(t+\tau))\big]\nonumber
\end{eqnarray}
Adding and subtracting the term $2\sum_{\tau=0}^{T_{\nu}-1}q_{S(t+\tau)}(\bv{U}(t))$ from the right hand side, we obtain:
\begin{eqnarray}
\|\bv{U}(t+T_{\nu})-\bv{U}^*_V\|^2\leq \|\bv{U}(t)-\bv{U}^*_V\|^2 + 2T_{\nu}B^2\label{eq:Tslotdrift_pre_manipulate}\qquad\quad\\
-2\sum_{\tau=0}^{T_{\nu}-1}\big[q_{S(t+\tau)}(\bv{U}^*_V)-q_{S(t+\tau)}(\bv{U}(t))\big]\qquad\nonumber\\
+2\sum_{\tau=0}^{T_{\nu}-1}\big[q_{S(t+\tau)}(\bv{U}(t+\tau))-q_{S(t+\tau)}(\bv{U}(t))\big].\nonumber
\end{eqnarray}


Since $\|\bv{U}(t)-\bv{U}(t+\tau)\|\leq\tau B$ and $\|\bv{A}^{(s_i)}(t)-\bv{\mu}^{(s_i)}(t)\|\leq B$, using  (\ref{eq:qla_osm_relation_pre}) and the fact that for any two vectors $\bv{x}$ and $\bv{y}$, $\bv{x}^T\bv{y}\leq\|\bv{x}\|\|\bv{y}\|$, we have:
\begin{eqnarray}
q_{S(t+\tau)}(\bv{U}(t+\tau))-q_{S(t+\tau)}(\bv{U}(t))\leq \tau B^2.\label{eq:qsi_boundedslope}
\end{eqnarray}
Hence:
\begin{eqnarray*}
\sum_{\tau=0}^{T_{\nu}-1}\big[q_{S(t+\tau)}(\bv{U}(t+\tau))-q_{S(t+\tau)}(\bv{U}(t))\big]\qquad\qquad\\
\leq\sum_{\tau=0}^{T_{\nu}-1}\big(\tau B^2\big)=\frac{1}{2}(T_{\nu}^2B^2-T_{\nu}B^2).
\end{eqnarray*}
Plug this into (\ref{eq:Tslotdrift_pre_manipulate}), we have:
\begin{eqnarray}
&&\hspace{-.4in}\|\bv{U}(t+T_{\nu})-\bv{U}^*_V\|^2 \leq \|\bv{U}(t)-\bv{U}^*_V\|^2 + (T^2_{\nu}+T_{\nu})B^2\label{eq:TMslotdrift_main}\\
&&\qquad\qquad-2\sum_{\tau=0}^{T_{\nu}-1}\big[q_{S(t+\tau)}(\bv{U}^*_V)-q_{S(t+\tau)}(\bv{U}(t))\big].\nonumber
\end{eqnarray}
Now denote $\script{Z}(t)=(\script{H}(t), \bv{U}(t))$, i.e., the pair of the history up to time $t$, $\script{H}(t)=\{S(\tau)\}^{t-1}_{\tau=0}$ and the current backlog. Taking expectations on both sides of (\ref{eq:TMslotdrift_main}), conditioning on $\script{Z}(t)$, we have:
\begin{eqnarray*}
&&\hspace{-.25in}\expect{\|\bv{U}(t+T_{\nu})-\bv{U}^*_V\|^2\left |\right. \script{Z}(t)}\qquad\qquad\qquad\quad \\
&&\qquad\leq \expect{\|\bv{U}(t)-\bv{U}^*_V\|^2 \left|\right. \script{Z}(t)} + (T^2_{\nu}+T_{\nu})B^2\label{eq:TMslotdrift_expect}\quad\\
&&\qquad-2\expect{\sum_{\tau=0}^{T_{\nu}-1}\big[q_{S(t+\tau)}(\bv{U}^*_V) 
-q_{S(t+\tau)}(\bv{U}(t))\big]\left|\right. \script{Z}(t)}.\nonumber
\end{eqnarray*}

Since the number of times $q_{s_i}(\bv{U})$ appears in the interval $[t, t+T_{\nu}-1]$ is $\|\script{T}_{s_i}(t, T_{\nu})\|$, we can rewrite the above as:
\begin{eqnarray}
&&\hspace{-.3in}\expect{\|\bv{U}(t+T_{\nu})-\bv{U}^*_V\|^2\left |\right. \script{Z}(t)}\qquad\qquad\qquad\quad\nonumber\\
&&\hspace{-.3in}\qquad\leq \expect{\|\bv{U}(t)-\bv{U}^*_V\|^2 \left|\right. \script{Z}(t)} + (T^2_{\nu}+T_{\nu})B^2\quad\nonumber\\
&&\hspace{-.3in}\qquad-2T_{\nu}\expect{\sum_{i=1}^{M} \frac{\|\script{T}_{s_i}(t, T_{\nu})\|}{T_{\nu}}\big[q_{s_i}(\bv{U}^*_V)-q_{s_i}(\bv{U}(t))\big]\left|\right. \script{Z}(t)}.\nonumber
\end{eqnarray}
Adding and subtracting $2T_{\nu}\sum_{i=1}^{M}p_{s_i}\big[q_{s_i}(\bv{U}^*_V)-q_{s_i}(\bv{U}(t))\big]$ from the right hand side, we have:
\begin{eqnarray}
&&\hspace{-.6in}\expect{\|\bv{U}(t+T_{\nu})-\bv{U}^*_V\|^2\left |\right. \script{Z}(t)}\label{eq:Tslotdrift_pre_almostdone}\\
&&\hspace{-.3in}\leq \expect{\|\bv{U}(t)-\bv{U}^*_V\|^2 \left|\right.  \script{Z}(t)} + (T^2_{\nu}+T_{\nu})B^2\quad\,\nonumber\\
&&\hspace{-.3in}\quad-2T_{\nu}\sum_{i=1}^{M}p_{s_i}\big[q_{s_i}(\bv{U}^*_V)-q_{s_i}(\bv{U}(t))\big]\quad\quad\nonumber\\
&&\hspace{-.3in}\quad-2T_{\nu}\expect{\sum_{i=1}^{M} \bigg[\frac{\|\script{T}_{s_i}(t, T_{\nu})\|}{T_{\nu}}-p_{s_i} \bigg]\times\quad\quad\nonumber\\
&&\qquad\qquad\big[q_{s_i}(\bv{U}^*_V)-q_{s_i}(\bv{U}(t))\big]\left|\right. \script{Z}(t)}.\nonumber
\end{eqnarray}
Denote the term inside the last expectation of (\ref{eq:Tslotdrift_pre_almostdone}) as $\script{Q}$, i.e., 
\begin{eqnarray}
\script{Q}=\sum_{i=1}^{M} \bigg[\frac{\|\script{T}_{s_i}(t, T_{\nu})\|}{T_{\nu}}-p_{s_i} \bigg]\big[q_{s_i}(\bv{U}^*_V)-q_{s_i}(\bv{U}(t))\big].
\end{eqnarray}
Using the fact that $q_{s_i}(\bv{U}^*_V)-q_{s_i}(\bv{U}(t))$ is a constant given $\script{Z}(t)$, we have:
\begin{eqnarray}
&&\hspace{-.3in}\expect{\script{Q}\left|\right. \script{Z}(t)}\nonumber\\
&&=\sum_{i=1}^{M} \bigg[\frac{\expect{\|\script{T}_{s_i}(t, T_{\nu})\|\left|\right.\script{Z}(t)}}{T_{\nu}}-p_{s_i} \bigg]\qquad\qquad\nonumber\\
&&\qquad\qquad\qquad\qquad\times\big[q_{s_i}(\bv{U}^*_V)-q_{s_i}(\bv{U}(t))\big]\nonumber\\
&&\leq \sum_{i=1}^{M} \bigg|\frac{\expect{\|\script{T}_{s_i}(t, T_{\nu})\|\left|\right.\script{Z}(t)}}{T_{\nu}}-p_{s_i} \bigg|\qquad\qquad\nonumber\\
&&\qquad\qquad\qquad\qquad\times\big|q_{s_i}(\bv{U}^*_V)-q_{s_i}(\bv{U}(t))\big|\nonumber
\end{eqnarray}
By (\ref{eq:qla_osm_relation_pre})$, q_{s_i}(\bv{U}^*_V)-q_{s_i}(\bv{U}(t))\leq B\|\bv{U}^*_V-\bv{U}(t)\|$, thus we have:
\begin{eqnarray}
\hspace{-.3in}\expect{\script{Q}\left|\right. \script{Z}(t)} &\leq&  B\|\bv{U}^*_V-\bv{U}(t)\| \nonumber\\
&&\times\sum_{i=1}^{M} \bigg|\frac{\expect{\|\script{T}_{s_i}(t, T_{\nu})\|\left|\right.\script{Z}(t)}}{T_{\nu}}-p_{s_i} \bigg|\nonumber\\
&\leq& \nu B\|\bv{U}^*_V-\bv{U}(t)\|,\label{eq:Tslotdrift_pre_almostdone_theexpect}
\end{eqnarray}
where the last step follows from the definition of $T_{\nu}$. Now by (\ref{eq:dual_separable}) and (\ref{eq:dual_singlestate}): \[\sum_{i=1}^{M}p_{s_i}\big[q_{s_i}(\bv{U}^*_V)-q_{s_i}(\bv{U}(t))\big]=q(\bv{U}^*_V)-q(\bv{U}(t)).\] 
Plug this and (\ref{eq:Tslotdrift_pre_almostdone_theexpect}) into (\ref{eq:Tslotdrift_pre_almostdone}),we have:
\begin{eqnarray*}
&&\hspace{-.6in}\expect{\|\bv{U}(t+T_{\nu})-\bv{U}^*_V\|^2\left |\right.  \script{Z}(t)}\qquad\qquad\qquad \\
&&\hspace{-.3in}\leq \expect{\|\bv{U}(t)-\bv{U}^*_V\|^2 \left|\right.  \script{Z}(t)} + (T^2_{\nu}+T_{\nu})B^2\qquad\\
&&\hspace{-.3in}-2T_{\nu}\big(q(\bv{U}^*_V)-q(\bv{U}(t))\big)+2T_{\nu} \nu B\|\bv{U}^*_V-\bv{U}(t)\|\nonumber
\end{eqnarray*}
Recall that $\script{Z}(t)=(\script{H}(t),\bv{U}(t))$. Taking expectation over $\script{H}(t)$ on both sides proves the lemma.
%
\end{proof}

\section*{Appendix B -- Proof of (\ref{eq:drift_tildeY_combine})  }
Here we prove that for $\tilde{Y}(t)$ defined in the proof of part (b) of Theorem \ref{theorem:prob_multicon_polyhedral}, we have:
\begin{eqnarray*}
\Delta_{T_{\nu}}(\tilde{Y}(t)) &\leq&e^{2wT_{\nu}B}-\frac{w\eta}{2}e^{w\tilde{Y}(t)},
\end{eqnarray*}
for all $\tilde{Y}(t)\geq0$
\begin{proof}
If $\tilde{Y}(t)>T_{\nu}B $, denote $\delta(t)=\tilde{Y}(t+T_{\nu})-\tilde{Y}(t)$. It is easy to see that $|\delta(t)|\leq T_{\nu}B$. Rewrite (\ref{eq:drift_exp_tildeY}) as:
\begin{eqnarray}
\Delta_{T_{\nu}}(\tilde{Y}(t))
&=&e^{w\tilde{Y}(t)}\expect{\big(e^{w\delta(t)}-1\big)\left|\right.\bv{U}(t)}.\label{eq:drift_exp_tildeY_delta}
\end{eqnarray}
By a Taylor expansion, we have that:
\begin{eqnarray}
e^{w\delta(t)}&=& 1+w\delta(t)+\frac{w^2\delta^2(t)}{2}g(w\delta(t)),\label{eq:e_delta_bound}
\end{eqnarray}
where $g(y)=2\sum_{k=2}^{\infty}\frac{y^{k-2}}{k!}=\frac{2(e^y-1-y)}{y^2}$ \cite{chung_concentration} has the following properties: 
\begin{enumerate}
\item $g(0)=1$; $g(y)\leq1$ for $y<0$; $g(y)$ is monotone increasing for $y\geq0$; 
\item For $y< 3$, 
\[g(y)=2\sum_{k=2}^{\infty}\frac{y^{k-2}}{k!}\leq \sum_{k=2}^{\infty}\frac{y^{k-2}}{3^{k-2}}=\frac{1}{1-y/3}.\]
\end{enumerate}
Thus by (\ref{eq:e_delta_bound}) we have: 
\begin{eqnarray}
e^{w\delta(t)}&\leq& 1+w\delta(t)+\frac{w^2T^2_{\nu}B^2}{2}g(wT_{\nu}B).\label{eq:e_delta_bound_g}
\end{eqnarray}
Plug this into (\ref{eq:drift_exp_tildeY_delta}), and note that $\tilde{Y}(t)>T_{\nu}B$, so  by (\ref{eq:Ytilde_sameasY}) we have $\expect{\delta(t)\left|\right.\bv{U}(t)}\leq-\eta$. Hence: 
\begin{eqnarray}
\hspace{-.2in}\Delta_{T_{\nu}}(\tilde{Y}(t)) &\leq&e^{w\tilde{Y}(t)}\big(-w\eta+\frac{w^2T^2_{\nu}B^2}{2}g(wT_{\nu}B)\big).\label{eq:drift_exp_tildeY_gB_alter}
\end{eqnarray}
Choosing $w=\frac{\eta}{T^2_{\nu}B^2+T_{\nu}B\eta/3}$, we see that $wT_{\nu}B<3$, thus:
\begin{eqnarray*}
\frac{w^2T^2_{\nu}B^2}{2}g(wT_{\nu}B)\leq \frac{w^2T^2_{\nu}B^2}{2}\frac{1}{1-wT_{\nu}B/3}=\frac{w\eta}{2},
\end{eqnarray*}
where the last equality follows since:
\begin{eqnarray*}
w=\frac{\eta}{T^2_{\nu}B^2+T_{\nu}B\eta/3} &\Rightarrow& w(T^2_{\nu}B^2+T_{\nu}B\eta/3)=\eta\\
&\Rightarrow& wT^2_{\nu}B^2=\eta-wT_{\nu}B\eta/3\\
&\Rightarrow& wT^2_{\nu}B^2\,\frac{1}{1-wT_{\nu}B/3}=\eta.
\end{eqnarray*}
Therefore (\ref{eq:drift_exp_tildeY_gB_alter}) becomes:
\begin{eqnarray}
\hspace{-.2in}\Delta_{T_{\nu}}(\tilde{Y}(t)) &\leq&-\frac{w\eta}{2}e^{w\tilde{Y}(t)}\leq e^{2wT_{\nu}B}-\frac{w\eta}{2}e^{w\tilde{Y}(t)}.
\end{eqnarray}
Now if $\tilde{Y}(t)\leq T_{\nu}B$, it is easy to see that $\Delta_{T_{\nu}}(\tilde{Y}(t))\leq e^{2wT_{\nu}B}-e^{w\tilde{Y}(t)}\leq e^{2wT_{\nu}B}-\frac{w\eta}{2}e^{w\tilde{Y}(t)}$, since $\tilde{Y}(t+T_{\nu})\leq T_{\nu}B+\tilde{Y}(t)\leq 2T_{\nu}B$ and $\frac{w\eta}{2}\leq1$, as $\eta<T_{\nu}B$. Therefore for all $\tilde{Y}(t)\geq0$, we see that (\ref{eq:drift_tildeY_combine}) holds.
\end{proof}

\section*{Appendix C-Proof of Lemma \ref{lemma:UtWt}}
Here we prove Lemma \ref{lemma:UtWt}.  To save space, we will sometimes use $[x]^+$ to denote $\max[x,0]$. 
\begin{proof}
It suffices to show that (\ref{eq:UWrelation}) holds for a single queue $j$. Also, when $\script{W}_j=0$, (\ref{eq:UWrelation}) trivially holds, thus we only consider $\script{W}_j>0$.

\emph{Part (A):} We first prove $U_j(t)\leq\max[W_j(t)-\script{W}_j, 0]+\delta_{max}$. First we see that it holds at $t=0$, since $W_j(0)=\script{W}_j$ and $U_j(t)=0$. It also holds for $t=1$. Since $U_j(0)=0$ and $W_j(0)=\script{W}_j$, we have $U_j(1)=A_j(0)\leq\delta_{max}$. Thus we have $U_j(1)\leq\max[W_j(1)-\script{W}_j, 0]+\delta_{max}$. 

Now assume $U_j(t)\leq\max[W_j(t)-\script{W}_j, 0]+\delta_{max}$ holds for $t=0,1,2,..., k$, we want to show that it also holds for $t=k+1$.
We first note that if $U_j(k)\leq\mu_j(k)$, the the result holds since then $U_j(k+1)=[U_j(k)-\mu_j(k)]^++A_j(k) =A_j(k)\leq\delta_{max}$. Thus we will consider $U_j(k)\geq\mu_j(k)$ in the following:

(A-I) Suppose $W_j(k)\geq\script{W}_j$. Note that in this case we have: 
\begin{eqnarray}
U_j(k)\leq W_j(k)-\script{W}_j+\delta_{max}.\label{eq:Ukknown1}
\end{eqnarray} 
Also, $U_j(t+1)=\max[U_j(t)-\mu_j(t),0]+A_j(t)$. Since $U_j(k)\geq\mu_j(k)$, we have:
\begin{eqnarray*}
\hspace{-.1in}U_j(k+1)&=&U_j(k)-\mu_j(k)+A_j(k)\\
\hspace{-.1in}&\leq& W_j(k)-\script{W}_j+\delta_{max}-\mu_j(k)+A_j(k)\\
\hspace{-.1in}&\leq& [W_j(k)-\mu_j(k)+A_j(k)-\script{W}_j]^++\delta_{max}\\
\hspace{-.1in}&\leq& \big[[W_j(k)-\mu_j(k)]^++A_j(k)-\script{W}_j\big]^++\delta_{max}\\
\hspace{-.1in}&=&\max[W_j(k+1)-\script{W}_j, 0]+\delta_{max},
\end{eqnarray*}
where the first inequality is due to (\ref{eq:Ukknown1}), the second and third inequalities are due to the $[x]^+$ operator, and the last equality follows from the definition of $W_j(k+1)$.


(A-II) Now suppose $W_j(k)<\script{W}_j$. In this case we have $U_j(k)\leq \delta_{max}$, $\tilde{A}_j(k)=[A_j(k)-\script{W}_j+W_j(k)]^+$ and:
\begin{eqnarray*}
U_j(k+1)=[U_j(k)-\mu_j(k)]^++\tilde{A}_j(k).
\end{eqnarray*}

First consider the case when $W_j(k)<\script{W}_j-A_j(k)$. In this case $\tilde{A}_j(k)=0$, so we have: 
\begin{eqnarray*}
U_j(k+1)=U_j(k)-\mu_j(k)\leq\delta_{max}-\mu_j(k)\leq\delta_{max},
\end{eqnarray*}
which implies $U_j(k+1)\leq\max[W_j(k+1)-\script{W}_j, 0]+\delta_{max}$. Else if $\script{W}_j-A_j(k)\leq W_j(k)< \script{W}_j$, we have:
\begin{eqnarray*}
U_j(k+1)&=&U_j(k)-\mu_j(k)+A_j(k)-\script{W}_j+W_j(k)\\
&\leq& W_j(k)-\script{W}_j+\delta_{max}-\mu_j(k)+A_j(k)\\
&\leq&\max[W_j(k+1)-\script{W}_j, 0]+\delta_{max},
\end{eqnarray*}
where the first inequality uses $U_j(k)\leq\delta_{max}$ and the second inequality follows as in (A-I). 


\emph{Part (B):} We now show that $U_j(t)\geq\max[W_j(t)-\script{W}_j, 0]$. First we see that it holds for $t=0$ since $W_j(0)=\script{W}_j$. We also have for $t=1$ that:
\begin{eqnarray*}
[W_j(1)-\script{W}_j]^+&=&\big[[W_j(0)-\mu_j(0)]^++A_j(0)-\script{W}_j\big]^+\\
&\leq& \big[[W_j(0)-\mu_j(0)-\script{W}_j]^++A_j(0)\big]^+\\
&=& A_j(0)
\end{eqnarray*}
Thus $U_j(1)\geq\max[W_j(1)-\script{W}_j, 0]$ since $U_j(1)=A_j(0)$. Now suppose $U_j(t)\geq\max[W_j(t)-\script{W}_j, 0]$ holds for $t=0,1,...,k$, we will show that it holds for $t=k+1$. We note that if $W_j(k+1)<\script{W}_j$, then $\max[W_j(k+1)-\script{W}_j, 0]=0$ and we are done. So we consider $W_j(k+1)\geq\script{W}_j$.

(B-I) First if $W_j(k)\geq\script{W}_j$, we have $\tilde{A}_j(k)=A_j(k)$. Hence:
\begin{eqnarray*}
[W_j(k+1)-\script{W}_j]^+
&=&[W_j(k)-\mu_j(k)]^++A_j(k)-\script{W}_j\\
&\leq& [W_j(k)-\mu_j(k)-\script{W}_j]^++A_j(k)\\
&\leq& [[W_j(k)-\script{W}_j]^+-\mu_j(k)]^++A_j(k)\\
&\leq& [U_j(k)-\mu_j(k)]^++A_j(k),
\end{eqnarray*}
where the first two inequalities are due to the $[x]^+$ operator and the last inequality is due to $U_j(k)\geq[W_j(k)-\script{W}_j]^+$. This implies $[W_j(k+1)-\script{W}_j]^+\leq U_j(k+1)$.

(B-II) Suppose $W_j(k)<\script{W}_j$. Since $W_j(k+1)\geq\script{W}_j$, we have $\script{W}_j-A_j(k)\leq W_j(k)< \script{W}_j$, for otherwise $W_j(k)<\script{W}_j-A_j(k)$ and $W_j(k+1)=[W_j(k)-\mu_j(t)]^++A_j(t)<\script{W}_j$. Hence in this case $\tilde{A}_j(k)=A_j(k)-\script{W}_j+W_j(k)\geq0$. 
\begin{eqnarray*}
&&[W_j(k+1)-\script{W}_j]^+\\
&=&[W_j(k)-\mu_j(k)]^++A_j(k)-\script{W}_j\\
&\leq&[W_j(k)+U_j(k)-\mu_j(k)]^++A_j(k)-\script{W}_j\\
&\leq&[U_j(k)-\mu_j(k)]^++A_j(k)-\script{W}_j+W_j(k)\\
&=& U_j(k+1)
\end{eqnarray*}
where the two inequalities are due to the fact that $U_j(k)\geq0$ and $W_j(k)\geq0$. 
\end{proof}

\section*{Appendix D-Proof of Lemma \ref{lemma:qla_tendtoopt}}
Here we prove Lemma \ref{lemma:qla_tendtoopt}. Recall that we use $\bv{x}_U$ to denote the vector $(x_{U}^{(s_1)}, x_{U}^{(s_2)}, ..., x_{U}^{(s_M)})^T$ chosen by OSM for a given $U(t)$, i.e., $\bv{x}_{U}$ achieves the infimum of (\ref{eq:dual}) at $U(t)$.
\begin{proof}
Now from the definition of $q(U(t))$, we have: 
\begin{eqnarray}
q(U(t))&=&\script{F}(\bv{x}_U)+U(t)\big[\script{G}_1(\bv{x}_U)-\script{B}_1(\bv{x}_U)\big]\nonumber\\
&=&\script{F}(\bv{x}_U)+U^*_V\big[\script{G}_1(\bv{x}_U)-\script{B}_1(\bv{x}_U)\big]\label{eq:tendtoopt_step}\\
&&\qquad\quad+(U(t)-U_V^*)\big[\script{G}_1(\bv{x}_U)-\script{B}_1(\bv{x}_U)\big].\nonumber
\end{eqnarray}
Using the fact that $q(U(t))<q(U^*_V)$ for $U(t)\neq U^*_V$, we have:
\begin{eqnarray}
q(U^*_V)&>&\script{F}(\bv{x}_U)+U^*_V\big[\script{G}_1(\bv{x}_U)-\script{B}_1(\bv{x}_U)\big]\label{eq:tendtoopt_step}\\
&&\qquad\quad+(U(t)-U_V^*)\big[\script{G}_1(\bv{x}_U)-\script{B}_1(\bv{x}_U)\big].\nonumber
\end{eqnarray}
This then implies:
\begin{eqnarray}
(U(t)-U_V^*)\big[\script{G}_1(\bv{x}_U)-\script{B}_1(\bv{x}_U)\big]\label{eq:tendtoopt_dual}\qquad\qquad\qquad\qquad\\
 <q(U^*_V)-\script{F}(\bv{x}_U)-U^*_V\big[\script{G}_1(\bv{x}_U)-\script{B}_1(\bv{x}_U)\big].\nonumber
\end{eqnarray}
However,  since: 
\[q(U^*_V)=\inf_{x^{(s_i)}\in\script{X}^{(s_i)}}\big\{\script{F}(\bv{x})+U^*_V\big[\script{G}_1(\bv{x})-\script{B}_1(\bv{x})\big]\big\},\] 
we have the right hand side of (\ref{eq:tendtoopt_dual}) being non-positive. Therefore:
\begin{eqnarray}
(U(t)-U_V^*)\big[\script{G}_1(\bv{x}_U)-\script{B}_1(\bv{x}_U)\big]<0.\label{eq:tendtoopt_dual_neg}
\end{eqnarray}
This proves (b). Now note that  under QLA, if the network state is $s_i$ then the chosen action $x^{(s_i)}_{U}$ minimizes: 
\begin{eqnarray}
Vf(s_i, x^{(s_i)})+U(t)\big[g_1(s_i, x^{(s_i)})-b_1(s_i, x^{(s_i)})\big],\label{eq:qla_min_isub}
\end{eqnarray}
over $\script{X}^{(s_i)}$ for the given $U(t)$. Therefore given $U(t)$, the expected value of the above quantity, i.e.,
\[\sum_ip_{s_i}\bigg\{Vf(s_i, x^{(s_i)})+U(t)\big[g_1(s_i, x^{(s_i)})-b_1(s_i, x^{(s_i)})\big]\bigg\},\] 
is minimized under QLA. Compare this fact to the definition of $q(U)$ in (\ref{eq:dual_separable}), we see that under QLA:
\begin{eqnarray}
q(U(t))=\expect{Vf(s_i, x^{(s_i)})\label{eq:dual_expectation_relation}\qquad\qquad\qquad\qquad\qquad\qquad\\
+U(t)\big[g_1(s_i, x^{(s_i)})-b_1(s_i, x^{(s_i)})\big]\left|\right.U(t)}.\nonumber
\end{eqnarray}
Thus similar as (\ref{eq:tendtoopt_step}), we have:
\begin{eqnarray}
\hspace{-.3in}&&q(U^*_V)>\expect{Vf(s_i, x^{(s_i)})\label{eq:dual_expection_relation_qla}\\
\hspace{-.3in}&&\qquad\qquad+U^*_V\big[g_1(s_i, x^{(s_i)})-b_1(s_i, x^{(s_i)})\big]\left|\right.U(t)}\nonumber\\
\hspace{-.3in}&&\qquad\qquad+(U(t)-U_V^*)\expect{g_1(s_i, x^{(s_i)})-b_1(s_i, x^{(s_i)})\left|\right.U(t)}.\nonumber
\end{eqnarray}
Now by (\ref{eq:dual_expectation_relation}) we see that $q(U^*_V)$ is the minimum of the expected value of (\ref{eq:qla_min_isub}) given $U^*_V$, we have: 
\begin{eqnarray}
q(U^*_V)&\leq& \expect{Vf(s_i, x^{(s_i)})\label{eq:dual_exp_less_qla}\\
&&+U^*_V\big[g_1(s_i, x^{(s_i)})-b_1(s_i, x^{(s_i)})\big]\left|\right.U(t)}.\nonumber
\end{eqnarray}
Subtract the right hand side of (\ref{eq:dual_exp_less_qla}) from both sides of (\ref{eq:dual_expection_relation_qla}) and use (\ref{eq:dual_exp_less_qla}), we see that Part (a) follows.
\end{proof}


\section*{Appendix E--Proof of Lemma \ref{lemma:moretendtoopt} and \ref{lemma:subgradient_qla}}
\begin{proof} (Lemma \ref{lemma:moretendtoopt})
We will prove the case when $0\leq U_1<U_2<U_V^*$, the other case can be similarly proven. First we have the following for the dual function:
\begin{eqnarray}
q(U_1)&=&\script{F}(\bv{x}_{U_1})+U_1\big[\script{G}_1(\bv{x}_{U_1})-\script{B}_1(\bv{x}_{U_1})\big]\label{eq:largerdrift_proof1}\\
&=&\script{F}(\bv{x}_{U_1})+U_2\big[\script{G}_1(\bv{x}_{U_1})-\script{B}_1(\bv{x}_{U_1})\big]\nonumber\\
&&\qquad\qquad\qquad+(U_1-U_2)\big[\script{G}_1(\bv{x}_{U_1})-\script{B}_1(\bv{x}_{U_1})\big].\nonumber
\end{eqnarray}
From the definition of $q(U_2)$ and $\bv{x}_{U_2}$, we see that:
\begin{eqnarray}
q(U_2)&=&\script{F}(\bv{x}_{U_2})+U_2\big[\script{G}_1(\bv{x}_{U_2})-\script{B}_1(\bv{x}_{U_2})\big]\nonumber\\
&\leq&\script{F}(\bv{x}_{U_1})+U_2\big[\script{G}_1(\bv{x}_{U_1})-\script{B}_1(\bv{x}_{U_1})\big].\label{eq:largerdrift_proof2}
\end{eqnarray}
Plug (\ref{eq:largerdrift_proof2}) into (\ref{eq:largerdrift_proof1}), we have: 
\begin{eqnarray}
q(U_1)&\geq&\script{F}(\bv{x}_{U_2})+U_2\big[\script{G}_1(\bv{x}_{U_2})-\script{B}_1(\bv{x}_{U_2})\big]\nonumber\\
&&\qquad\qquad\qquad+(U_1-U_2)\big[\script{G}_1(\bv{x}_{U_1})-\script{B}_1(\bv{x}_{U_1})\big]\nonumber\\
&=&\script{F}(\bv{x}_{U_2})+U_1\big[\script{G}_1(\bv{x}_{U_2})-\script{B}_1(\bv{x}_{U_2})\big]\nonumber\\
&&+(U_1-U_2)\bigg\{\big[\script{G}_1(\bv{x}_{U_1})-\script{B}_1(\bv{x}_{U_1})\big]\label{eq:largerdrift_proof3}\\
&&\qquad\qquad\qquad\qquad-\big[\script{G}_1(\bv{x}_{U_2})-\script{B}_1(\bv{x}_{U_2})\big]\bigg\}.\nonumber
\end{eqnarray}
Now similar as in (\ref{eq:largerdrift_proof2}) we have $q(U_1)\leq \script{F}(\bv{x}_{U_2})+U_1\big[\script{G}_1(\bv{x}_{U_2})-\script{B}_1(\bv{x}_{U_2})\big]$.  Therefore from (\ref{eq:largerdrift_proof3}) we obtain:
\begin{eqnarray*}
0\geq(U_1-U_2)\bigg\{\big[\script{G}_1(\bv{x}_{U_1})-\script{B}_1(\bv{x}_{U_1})\big]\qquad\qquad\qquad\qquad\\
-\big[\script{G}_1(\bv{x}_{U_2})-\script{B}_1(\bv{x}_{U_2})\big]\bigg\}.
\end{eqnarray*}
Since $U_1<U_2$, $\script{G}_1(\bv{x}_{U_1})-\script{B}_1(\bv{x}_{U_1})\geq\script{G}_1(\bv{x}_{U_2})-\script{B}_1(\bv{x}_{U_2})$. Similar as in the proof of Lemma \ref{lemma:qla_tendtoopt}, we see that we also have:
\begin{eqnarray*}
\expect{g_{1}(s_i, x^{(s_i)}_{U_1})-b_{1}(s_i, x^{(s_i)}_{U_1})\left|\right.U_1}\qquad\qquad\qquad\qquad\\
\geq \expect{g_{1}(s_i, x^{(s_i)}_{U_2})-b_{1}(s_i, x^{(s_i)}_{U_2})\left|\right.U_2}.
\end{eqnarray*}
From Lemma \ref{lemma:qla_tendtoopt} Part (a) we see that they are both positive. 
\end{proof}

\begin{proof} (Lemma \ref{lemma:subgradient_qla}) 
Note that from (\ref{eq:qla_osm_relation_pre}), we have:
\begin{eqnarray*}
(\bv{U}_V^*-\bv{U}(t))^{T}\big[\sum_ip_{s_i}(\bv{A}^{(s_i)}(t)-\bv{\mu}^{(s_i)}(t))\big]\qquad\qquad\\
\geq q(\bv{U^*_V})-q(\bv{U}(t)).
\end{eqnarray*}
This leads to the following inequality:
\begin{eqnarray*}
\sum_{j=1}^r(U_{Vj}^*-U_j(t))\expect{\big[g_j(s_i, x^{(s_i)})-b_j(s_i, x^{(s_i)})\big]\left|\right. \bv{U}(t)}\\
\geq q(\bv{U^*_V})-q(\bv{U}(t)).\nonumber
\end{eqnarray*}
Taking $r=1$, we see that Lemma \ref{lemma:subgradient_qla} follows.
\end{proof}

$\vspace{-.3in}$
\bibliographystyle{unsrt}
\bibliography{../mybib}

\end{document}